\documentclass[a4paper]{amsart} %%%
\usepackage[margin=1in]{geometry}	%%%%	\thispagestyle{empty}	Luba:[margin=3cm]			%%%%
\usepackage{amsmath}	%%%%%%%%%%%%%%%%	load before amsthm									%%%%
\usepackage{amsthm}	%%%%%%%%%%%%%%%%%%%%	load after amsmath									%%%%
\usepackage{amssymb}	%%%%%%%%%%%%%%%%														%%%%
\usepackage{extarrows}	%%%%%%%%%%%%%%%%	\x...arrow[subtext]{supertext}						%%%%
\usepackage{mathtools}	%%%%%%%%%%%%%%%%														%%%%
\usepackage{thmtools}	%%%%%%%%%%%%%%%%														%%%%
\usepackage{amsfonts}	%%%%%%%%%%%%%%%%														%%%%
\usepackage[scr=boondoxo]{mathalfa}	%%%%	\mathscr uses boondoxo font							%%%%
\usepackage{fontawesome}	%%%%%%%%%%%%														%%%%
\usepackage{bm}	%%%%%%%%%%%%%%%%%%%%%%%%	load after other math fonts packages				%%%%
\usepackage[shortlabels]{enumitem}	%%%%	[before=\leavevmode] for start of theorems			%%%%
\usepackage{graphicx}	%%%%%%%%%%%%%%%%														%%%%
\usepackage[font={small}]{caption}	%%%%	\captionsetup{labelformat=empty,labelsep=none}		%%%%
\usepackage{subcaption}	%%%%%%%%%%%%%%%%	\captionsetup[figure]{labelfont=it}					%%%%
\usepackage[dvipsnames]{xcolor}	%%%%%%%%	load before tikz									%%%%
\usepackage{comment}	%%%%%%%%%%%%%%%%														%%%%
\usepackage{etoolbox}	%%%%%%%%%%%%%%%%														%%%%
\usepackage{lmodern}	%%%%%%%%%%%%%%%%														%%%%
\usepackage{upgreek}	%%%%%%%%%%%%%%%%														%%%%
\usepackage[doi=true,isbn=true,giveninits=true,maxbibnames=99,backend=biber,url=true]{biblatex}	%%%%
\usepackage[hypertexnames=false]{hyperref}	%%%%%%%%%%%%%%%%	[bookmarksdepth=4]									%%%%
\usepackage{xurl}	%%%%%%%%%%%%%%%%%%%%														%%%%
\usepackage[capitalise]{cleveref}	%%%%	Luba:[nameinlink]									%%%%
\hypersetup{colorlinks=true,linktoc=section,urlcolor={Aquamarine!85!black},%					%%%%
	citecolor={YellowOrange!85!black},linkcolor={violet}}										%%
\crefname{thmabc}{Theorem}{Theorems}
\Crefname{thmabc}{Theorem}{Theorems}
\crefname{corabc}{Corollary}{Corollaries}
\Crefname{corabc}{Corollary}{Corollaries}
\makeatletter	%%%%%%%%%%%%%%%%%%%%%%%%	cleveref compatibility fix		--DO NOT REMOVE--	%%%%
\AtBeginDocument{%					%%%%									--DO NOT REMOVE--	%%
\patchcmd\label@noarg{\edef\@tempb}{\protected@edef\@tempb}{}{}}	%%		--DO NOT REMOVE--	%%
\makeatother	%%%%%%%%%%%%%%%%%%%%%%%%									--DO NOT REMOVE--	%%
\renewbibmacro*{in:}{%				%%%%	removes "In:" when journal entry empty				%%%%
	\iffieldundef{journaltitle}																	%%
	{}																							%%
	{\printtext{\bibstring{in}\intitlepunct}}													%%
}%																								%%
\theoremstyle{plain}	%%%%%%%%%%%%%%%%														%%%%
\newtheorem{thm}{Theorem}[section]	%%%%	[<counter>]=[chapter] when dc=book					%%%%
\crefname{thm}{Theorem}{Theorems}	%%%%	tells cleveref what is the singular and plural		%%
\newtheorem{lem}[thm]{Lemma}	%%%%%%%%	Luba:{}[<empty>]{}[section]							%%%%
\crefname{lem}{Lemma}{Lemmata}	%%%%%%%%	tells cleveref what is the singular and plural		%%
\newtheorem{prop}[thm]{Proposition}	%%%%														%%%%
\newtheorem{cor}[thm]{Corollary}	%%%%														%%%%
\theoremstyle{definition}	%%%%%%%%%%%%														%%%%
\crefname{eg}{Example}{Examples}	%%%%	tells cleveref what is the singular and plural		%%
\theoremstyle{plain}	%%%%%%%%%%%%%%%%														%%%%
\newtheorem{rem}[thm]{Remark}	%%%%%%%%	Luba:{remark}										%%%%
\theoremstyle{plain}	%%%%%%%%%%%%%%%%														%%%%
\newtheorem*{thm*}{Theorem}	%%%%%%%%%%%%														%%%%
\newtheorem*{lem*}{Lemma}	%%%%%%%%%%%%														%%%%
\newtheorem*{prop*}{Proposition}	%%%%														%%%%
\newtheorem*{cor*}{Corollary}	%%%%%%%%														%%%%
\theoremstyle{definition}	%%%%%%%%%%%%														%%%%
\newtheorem*{dfn*}{Definition}	%%%%%%%%														%%%%
\newtheorem*{eg*}{Example}	%%%%%%%%%%%%	Luba:{exmp}											%%%%
\theoremstyle{remark}	%%%%%%%%%%%%%%%%														%%%%
\newtheorem*{rem*}{Remark}	%%%%%%%%%%%%														%%%%
\theoremstyle{plain}	%%%%%%%%%%%%%%%%														%%%%
\newtheorem{thmabc}{Theorem}	%%%%%%%%														%%%%
\newtheoremstyle{named}{}{}{\itshape}{}{\bfseries}{.}{.5em}{\thmnote{#3}}	%%					%%%%
\theoremstyle{named}	%%%%%%%%%%%%%%%%														%%
\numberwithin{equation}{section}	%%%%	<counter> resets on section							%%%%
\patchcmd{\subsection}{-.5em}{1sp}{}{}%%	skip to leave below heading							%%%%
\makeatletter	%%%%%%%%%%%%%%%%%%%%%%%%	\hrulefill											%%%%
\def\hrulefill{\noindent\leavevmode\leaders\hrule\hfill\kern\z@}	%%							%%
\makeatother	%%%%%%%%%%%%%%%%%%%%%%%%														%%
\DeclareMathOperator{\dist}{dist}	%%%%	distance											%%%%
\DeclareMathOperator{\diam}{diam}	%%%%	diameter		                                    %%%%

\DeclareMathOperator{\id}{id}	%%%%%%%%	identity operator									%%%%
\ExplSyntaxOn	%%%%%%%%%%%%%%%%%%%%%%%%	new mb.., mc.., ms..alphabets						%%%%
\NewDocumentCommand{\definealphabet}{mmmm}{%	#1 = prefix, #2 = command, #3 = start, #4 = end	%%
	\int_step_inline:nnn { `#3 } { `#4 }{%														%%
		\cs_new_protected:cpx { #1 \char_generate:nn { ##1 }{ 11 } }{%							%%
			\exp_not:N #2 { \char_generate:nn { ##1 } { 11 } }									%%
		}																						%%
	}																							%%
}																								%%
\ExplSyntaxOff																					%%
\definealphabet{mb}{\mathbb}{A}{Z}	%%%%														%%
\definealphabet{mc}{\mathcal}{A}{Z}	%%%%														%%
\definealphabet{ms}{\mathscr}{A}{Z}	%%%%														%%

\DeclareMathOperator{\jul}{\mathscr{J}}	%%														%%%%
\DeclareMathOperator{\fat}{\mathscr{F}}	%%														%%%%
\newcommand{\CC}{{\hat{\mathbb{C}}}}%%%%														%%%%
\newcommand{\seq}[1]{{\bm{{#1}}}}	%%%%														%%%%
\newcommand{\area}{m_2}

\title{Cantor Julia sets of positive area in non-autonomous dynamics}
\author[L.~Pathirana]{Lubashan Pathirana}
\address{Department of Mathematical Sciences, University of Copenhagen, Universitetsparken 5, 2100 Copenhagen \O, Denmark}
\email{lpk@math.ku.dk}
\begin{document}

%%%%%%%%%%%%%%%%%%%%%%%%%%%%%%%%%%%%%%%%%%%%%%%%%%%%%%%%%%%%%%%%%%%%%%%%%%%%%%%%%%%%%%%%%%%%%%%%%%%%
%%%%%%%%%%%%%%%%%%%%%%%%%%%%%%%%%%%%%%%%%%%%%%%%%%%%%%%%%%%%%%%%%%%%%%%%%%%%%%%%%%%%%%%%%%%%%%%%%%%%

\begin{abstract}
    For every \(a>0\), we prove the existence of a sequence of polynomials with rational coefficients whose non-autonomous iteration produces a Julia set of planar area \(a\).
    We also construct a sequence of rational maps with rational coefficients whose Julia set has the same prescribed area.
    No individual map in this rational sequence is M\"obius-conjugate to a polynomial, and no choice of M\"obius coordinates at successive phases turns it into a polynomial sequence.
    In both cases, for every \(n\geq1\), the Julia set obtained by starting the iteration at the \(n\)-th map, called the phase-\(n\) Julia set, is a Cantor set and meets every neighborhood of each of its points in positive planar area.
    For every \(A\in(0,1)\), we construct two sequences of quadratic polynomials, each with \(\infty\) in every phase Julia set and with first phase Julia set of normalized spherical area \(A\).
    In one sequence, every phase Julia set is a Cantor set of locally positive spherical area.
    In the other, every phase Julia set has countably infinitely many connected components and equals the closure of its interior.
    Both sequences admit invariant measurable line fields whose essential supports are the corresponding phase Julia sets.
    Finally, for every \(d\in[0,2]\), we construct a sequence of quadratic polynomials for which every phase Julia set is a Cantor set of zero spherical area with local Hausdorff dimension \(d\) at each of its points.
\end{abstract}

%%%%%%%%%%%%%%%%%%%%%%%%%%%%%%%%%%%%%%%%%%%%%%%%%%%%%%%%%%%%%%%%%%%%%%%%%%%%%%%%%%%%%%%%%%%%%%%%%%%%
%%%%%%%%%%%%%%%%%%%%%%%%%%%%%%%%%%%%%%%%%%%%%%%%%%%%%%%%%%%%%%%%%%%%%%%%%%%%%%%%%%%%%%%%%%%%%%%%%%%%

\maketitle
% \tableofcontents

%%%%%%%%%%%%%%%%%%%%%%%%%%%%%%%%%%%%%%%%%%%%%%%%%%%%%%%%%%%%%%%%%%%%%%%%%%%%%%%%%%%%%%%%%%%%%%%%%%%%
%%%%%%%%%%%%%%%%%%%%%%%%%%%%%%%%%%%%%%%%%%%%%%%%%%%%%%%%%%%%%%%%%%%%%%%%%%%%%%%%%%%%%%%%%%%%%%%%%%%%

\section{Introduction}
\label{sec:intro}

Given a rational map $f:\CC\to\CC$, where $\CC\:=\mbC\cup\{\infty\}$ is the Riemann sphere, with $\deg(f)\geq2$, autonomous rational dynamics studies repeated compositions $f^{\circ k}$ for $k\in\mbN$. 
The \emph{Fatou set} of the function $f$, denoted by $\fat(f)$,  is the set of points that have a neighborhood on which $\{f^{\circ k}:k\geq0\}$ is a normal family, with respect to the spherical metric.
Its complement $\jul(f):=\CC\setminus \fat(f)$ is the \emph{Julia set} of $f$. 
Both sets are completely invariant with respect to $f$: taking either the image or the preimage (in $\CC$) under $f$ leaves each set unchanged. 
Moreover, $\jul(f)$ is nonempty, compact, and perfect. It is also known that $\jul(f)$ is either the entire sphere $\CC$ or has no interior; see \cite{Beardon1991,Milnor2006}.

In the non-autonomous case, one starts with a sequence $\seq f = (f_n)_{n\ge 1}$ of nonconstant rational functions on $\CC$. We generally allow maps of degree one in this setting, and we shall state degree restrictions explicitly when needed.
Given any $n\in\mbN$, called the \emph{phase}, we consider the forward iterates
\[
    f_{n+k}\circ\cdots\circ f_n \qquad k\ge 0
\]
and we define 
\[
    \fat_n := \left\{z\in\CC : \left(f_{n+k}\circ\cdots\circ f_n\right)_{k\ge 0} \text{ is normal at } z\right\}
\]
as the Fatou set of $\seq{f}$ \emph{with phase $n$}, or simply the phase $n$ Fatou set of $\seq{f}$.
Here the normality is taken with respect to the spherical metric on $\CC$.
We use its equivalent local equicontinuity formulation and Montel's theorem for meromorphic families omitting three fixed values; see \cite{Milnor2006}.
Similar to the autonomous case, the Julia set with phase $n$ (or simply the phase $n$ Julia set) is defined as $\jul_n := \CC\setminus \fat_n$.
The collection $(\jul_n)_{n\ge 1}$ is referred to as the Julia sequence of $\seq{f}$ or as the iterated Julia sets of $\seq{f}$. 
The first phase Julia set $\jul_1$ is loosely referred to as \emph{the} Julia set of the sequence $\seq{f}$.

If infinitely many of the rational functions in the sequence $\seq f = (f_n)_{n\ge 1}$ have degree at least $2$, then the Julia sets $\jul_n$ for each phase $n$ are known to be non-empty. 
These sets $(\jul_n)_{n\ge 1}$ and  $(\fat_n)_{n\ge 1}$ satisfy the following invariance property; 
See \cite[Theorem~3.1, Proposition~3.3]{monomial}; see also \cite[Theorem~1.1]{Comerford2013}.
Let $\seq{f}=(f_1,f_2,\dots)$ be a sequence of nonconstant rational maps of $\CC$.
Then for every~$n\ge 1$
\begin{equation}
\label{thm:invariance deterministic}
  \left(f_n \circ \cdots \circ f_1\right)\left(\jul_1\right) = \jul_{n+1}
    \qquad\text{and}\qquad
    \left(f_n \circ \cdots \circ f_1\right)^{-1} (\jul_{n+1}) =\jul_1.
\end{equation}
The same equalities hold for the Fatou sets.
Furthermore, if $\seq f$ contains infinitely many rational functions of degree at least $2$, then for each $n$ the Julia set $\jul_n$ at phase $n$ is non-empty.

In autonomous dynamics, a quadratic polynomial can have a Julia set of positive planar area and empty interior~\cite{Avila_2022,Buff2012}, and the Julia set of a transcendental entire function can contain a Cantor set of positive area~\cite[Example~4]{Eremenko_1987}.
In non-autonomous dynamics, a bounded quadratic sequence can have Julia sets of positive area at every phase~\cite[Theorem~1.2]{Comerford2013}, and positive area occurs even for maps $z^2+c_n$ with $c_n\in\{0,\pm1/4,1/2\}$~\cite[p.~361]{Bruck01}.
It follows from \cite[Section~5]{monomial} that for every $\alpha\in[0,+\infty]$ there is a monomial sequence with $\area(\jul_1)=\alpha$ and $\operatorname{int}(\jul_1)=\emptyset$; moreover (while it is not stated as such in \cite{monomial}), the monomial sequence can be taken so that it is quadratic and each monomial in the sequence has rational coefficients.  
Related work addresses the dimension of autonomous Cantor Julia sets~\cite{Yang} and questions of disconnectedness, dimension, and measure for polynomial sequences~\cite{Alpan_2017,ComerfordStankewitzSumi2020,Lech2021}. The behavior of phase Julia sets under conjugacy and iteration has also been studied~\cite{selfsimilar,Comerford03}.
Motivated by these examples, we study polynomial and rational sequences whose phase Julia sets are Cantor sets of positive planar area.

\medskip 

For $z\in\mathbb{C}$ and $r>0$, let $B(z,r):=\{w\in\mathbb{C}: |w-z|<r\}$ denote the open disk of radius $r$ centered at $z$.
We write $\mbD := B(0,1)$ for the open unit disk, and we denote its closure by \(\overline{\mbD}\). 
We denote the closed disk of radius \(r\) centered at \(z\) by \(\overline{B(z,r)}\).
The Gaussian rational field \(\{u+iv:u,v\in\mbQ\}\) is denoted by \(\mbQ(i)\). 
For a field \(F\subseteq\mbC\), we write \(F[z]\) for the ring of polynomials with coefficients in \(F\), and \(F(z)\) for the field of quotients \(P(z)/Q(z)\), where \(P,Q\in F[z]\) and \(Q\) is not the
zero polynomial.
We denote by $\area$ the planar Lebesgue measure on $\CC$. 
In what comes below, the term \emph{Cantor set} refers to a nonempty, compact, perfect, and totally disconnected metrizable space.

\smallskip 

The first main contribution of this paper is the following result. 

\begin{thmabc}
\label{thm:main}
    \begin{enumerate}
        \item
            For every \(a>0\), there exists a sequence of polynomials \(\seq{f}=(f_n)_{n\geq1}\) in \(\mbQ[z]\) whose first phase Julia set satisfies \(\area(\jul_1)=a\).

        \item
            For every \(a>0\), there exists a sequence of rational maps \(\seq{f}=(f_n)_{n\geq1}\) in \(\mbQ(z)\), none of which  is M\"obius-conjugate to a polynomial, such that \(\area(\jul_1)=a\).
            Moreover, there is no sequence of M\"obius maps \((T_n)_{n\geq1}\) for which \(T_{n+1}\circ f_n\circ T_n^{-1}\) is a polynomial for every \(n\geq1\).
    \end{enumerate}
    In both cases, for every \(n\geq1\), the phase Julia set \(\jul_n\) is a Cantor subset of \(\mbC\) with positive planar area; in particular, it has empty interior.
    Moreover, every  neighborhood of each of its points meets it in positive planar area:
    \[
        \area\bigl(\jul_n\cap B(z,r)\bigr)>0
        \qquad\text{for every }z\in\jul_n\text{ and every }r>0.
    \]
\end{thmabc}

The proof of \Cref{thm:main} proceeds in several stages.

In \Cref{thm:polynomial} we first establish that, for every \(0<a<\pi\), there exists a sequence of polynomials in \(\mbQ[z]\) with \(\area(\jul_1)=a\) and with the stated Cantor and local area properties at every phase. Then \Cref{cor:polynomail} extends this result to obtain any area via a fixed dilation while preserving the desired properties required in part~(1) of \Cref{thm:main}

For part~(2), \Cref{thm:rational} gives a recursive formula for a sequence of rational maps in \(\mbQ(z)\) whose first phase Julia set has positive planar area.
This sequence has the required Cantor and local area properties and satisfies both assertions concerning M\"obius changes of coordinates. 
Without the restriction to rational coefficients, prescribing the area would follow immediately from a fixed dilation conjugacy. 
To prescribe the area exactly, with rational coefficients, \Cref{cor:prescribed-area} modifies the recursion by choosing a rational multiplier at each step to first obtain any area in $(0,\pi)$ and then scale it via an integer dilation.
Together, \Cref{cor:polynomail,cor:prescribed-area} prove \Cref{thm:main}.

Each of the two sequences in \Cref{cor:polynomail,cor:prescribed-area} satisfies the property that there is an open disk \(D\mathrel := B(0,R)\) such that
\[
    f_n^{-1}(\overline D)\Subset D
    \qquad\text{and}\qquad
    \jul_n\Subset D
\]
for each $n\in\mbN$.
Here \(A\Subset B\) means that \(\overline A\) is a compact subset of \(B\) and is referred to as $A$ is \emph{compactly contained} in $B$.
In \Cref{thm:polynomial,thm:rational} we may take \(D=\mbD\).
For sequences of polynomials of degree at least two, the condition above implies the common exterior escape condition \eqref{eq:uniform-exterior}. 
In this degree range, this common exterior escape condition is equivalent to the common invariant domain hypothesis used by B\"uger~\cite{selfsimilar} (see \Cref{rem:buger-common-escape}).

\medskip

Part~(1) of \Cref{thm:main} also holds outside the class of polynomial sequences satisfying \eqref{eq:uniform-exterior}. 
Starting with a sequence from \Cref{thm:polynomial}, \Cref{prop:exterior} applies a dilation at each phase. 
The resulting sequence has the same first phase Julia set, and each later phase Julia set is a dilation of the corresponding original set. For every \(R>0\), some map in the resulting sequence sends a point of \(\{|z|>R\}\) into \(\{|z|<R\}\). 
Applying the integer dilation from the proof of \Cref{cor:polynomail} to the sequence in \Cref{prop:exterior} gives, for every \(a>0\), a sequence in \(\mbQ[z]\) with \(\area(\jul_1)=a\) and with the local area property at every phase.

\medskip 

Having obtained the existence of a sequence of polynomials satisfying Part~(1) of \Cref{thm:main} above, we now seek an explicit sequence with the same Cantor and local positive-area properties. 
The polynomial sequences in \Cref{thm:polynomial,prop:exterior} have degrees tending to infinity, and \(\infty\) belongs to the Fatou set at every phase. 
In \Cref{sec:quadratic}, we give an explicit sequence of quadratic polynomials with \(\infty\) in every phase Julia set. 
To include neighborhoods of \(\infty\) in the area statements, we use spherical area \(m_{\mathrm{sph}}\) on \(\CC\), normalized by \(m_{\mathrm{sph}}(\CC)=1\).

\smallskip 

We now present our next result. 

\begin{thmabc}
\label{thm:main-quadratic}
    Fix the sequence of quadratic polynomials \(\seq{q} = (q_m)_{m\ge 1}\) in \(\mbQ[z]\) where 
    \[
        q_m(z) := z+2^{-2^{m^2+4}}z^2
        \qquad\text{ for each }\qquad m\in \mbN.
    \]
    Given a bounded sequence \((a_m)_{m\geq1}\) in \(\mbQ(i)\), put \(M_0:=\id\) and define
    \[
        M_m(z) := m(z-a_m),
        \qquad\text{ and }\qquad
        p_m := M_m\circ q_m\circ M_{m-1}^{-1}
        \qquad(m\geq1).
    \]
    Write \(\jul_n\) for the phase-\(n\) Julia set of \((p_m)_{m\geq1}\).
    The following choices of sequences  \((a_m)_{m\ge 1}\) are possible. 
    \begin{enumerate}
        \item
            For every \(A\in(0,1)\), there is a sequence \((a_m)_{m\ge 1}\) for which \(m_{\mathrm{sph}}(\jul_1)=A\) and every phase Julia set, \(\jul_n\), is a Cantor set of positive spherical area.
            In particular, at every phase \(n\), we have \(m_{\mathrm{sph}}(V\cap\jul_n)>0\) for every open set \(V\subseteq\CC\) meeting \(\jul_n\).

        \item
            For every \(A\in(0,1)\), there is a sequence \((a_m)_{m\ge 1}\) for which \(m_{\mathrm{sph}}(\jul_1)=A\) and every phase Julia set, \(\jul_n\), has countably infinitely many connected components.
            For every $n\in\mbN$, each open set \(V\subseteq\CC\) meeting \(\jul_n\) contains an open planar disk lying entirely in \(\jul_n\).
            Equivalently, \(\jul_n=\overline{\operatorname{int}(\jul_n)}\) at every phase \(n\in\mbN\), where the interior and closure are taken in \(\CC\) with spherical metric.

        \item
            For every \(d\in[0,2]\), there is a sequence \((a_m)_{m\ge 1}\) for which every phase $n$ Julia set, \(\jul_n\), is a Cantor set of zero spherical area and with local Hausdorff dimension $d$ at each of its points, i.e., for each $ n\in\mbN$, and for every open set \(V\subseteq\CC\) meeting \(\jul_n\), we have $\dim_{\mathrm H}(V\cap\jul_n)=d$.
    \end{enumerate}
    In all cases, \(p_m\in\mbQ(i)[z]\) is quadratic for every \(m\), and \(\infty\in\jul_n\) at every phase \(n\).
    In parts~(1) and~(2), the sequence \((p_m)_{m\geq1}\) also admits an invariant sequence of measurable line fields whose essential support at phase \(n\) is exactly \(\jul_n\).
\end{thmabc}

The three parts of \Cref{thm:main-quadratic} are obtained by applying the following general result in \Cref{thm:quadratic-compact-pullback} to different compact sets \(K\subset\mbC\).
The theorem identifies the finite part of each phase Julia set as the preimage of \(K\) under an entire function and relates its local area, Hausdorff dimension, and interior to those of \(K\).

Every nonempty compact set \(K\subseteq\mbC\) is the set of subsequential limits of a bounded sequence in \(\mbQ(i)\).
Indeed, for each \(j\geq1\), list the midpoints of all closed dyadic squares of side \(2^{-j}\) that meet \(K\).
Each list is finite and consists of points in \(\mbQ(i)\).
Every midpoint lies within \(\sqrt{2}\,2^{-j-1}\) of \(K\), and every point of \(K\) lies within the same distance of a midpoint in the list.
Concatenating these lists gives the desired sequence.
For any bounded sequence \((a_m)\) whose set of subsequential limits is \(K\), we also have \(\dist(a_m,K)\to0\), since otherwise a subsequence staying a positive distance from \(K\) would have a further subsequence converging outside \(K\).

\begin{restatable}{thmabc}{QuadraticCompactPullback}
\label{thm:quadratic-compact-pullback}
    Let \(K\subset\mbC\) be a nonempty compact set and choose any bounded sequence \((a_m)_{m\geq1}\) in \(\mbQ(i)\) whose set of subsequential limits is \(K\).
    Let \(q_m\), \(M_m\), \(p_m\), and \(\jul_n\) be as in \Cref{thm:main-quadratic}.

    At each phase \(n\), the compositions \(Q_{n,k} := q_{n+k-1}\circ\cdots\circ q_n\) converge locally uniformly on \(\mbC\), as \(k\to\infty\), to a nonconstant entire function \(U_n\).
    Put \(H_n := U_n\circ M_{n-1}^{-1}\).
    Then
    \[
        \jul_n=H_n^{-1}(K)\cup\{\infty\}.
    \]
    For every \(w\in\mbC\), the set \(H_n^{-1}(\{w\})\) is unbounded.

    For every open set \(V\subseteq\CC\) meeting \(\jul_n\), put \(W := H_n(V\cap\mbC)\).
    Then
    \[
        m_{\mathrm{sph}}(V\cap\jul_n)>0
        \qquad\text{if and only if}\qquad
        \area(W\cap K)>0.
    \]
    We also have \(\dim_{\mathrm H}(V\cap\jul_n)=\dim_{\mathrm H}(W\cap K)\).
    Moreover, a finite point \(z\in\mbC\) belongs to \(\operatorname{int}(\jul_n)\) if and only if \(H_n(z)\in\operatorname{int}(K)\).
\end{restatable}

In \Cref{sec:quadratic}, we first prove \Cref{thm:quadratic-compact-pullback} and then use it to establish \Cref{thm:main-quadratic}.
For part~(1), let \(\mcC\subset[0,1]\) be the Smith--Volterra--Cantor set.
Taking \(K := \mcC+i\mcC\) gives the explicit positive-area example in \Cref{prop:quadratic-fat}.
To prescribe the first phase spherical area \(A\), we instead choose a compact Cantor set \(K\) with \(m_{\mathrm{sph}}(U_1^{-1}(K))=A\), as in \Cref{cor:quadratic-prescribed-spherical-area}.
For part~(2), a family of disjoint closed disks accumulating at a point gives the component count in \Cref{rem:countable-julia-components}, and \Cref{cor:countable-prescribed-spherical-area} gives the prescribed first phase spherical area.
The invariant line fields in parts~(1) and~(2) follow from \Cref{cor:quadratic-line-fields}.
For part~(3), we choose a planar Cantor set \(K\) of zero area such that \(\dim_{\mathrm H}(W\cap K)=d\) for every open set \(W\subseteq\mbC\) meeting \(K\). The proof of part~(3) is carried out in \Cref{thm:quadratic-prescribed-dimension}.

Other choices of the compact set \(K\) in \Cref{thm:quadratic-compact-pullback} can yield further examples of Julia sets.
For example, \Cref{sec:quadratic-mixed} gives a quadratic sequence for which every phase Julia set has nonempty interior and contains two finite points \(x\) and \(y\) with the following properties.
Every sufficiently small disk about \(x\) meets the Julia set in a set of positive planar area and empty interior, while every sufficiently small disk about \(y\) meets it in a set of zero planar area and Hausdorff dimension one.

%%%%%%%%%%%%%%%%%%%%%%%%%%%%%%%%%%%%%%%%%%%%%%%%%%%%%%%%%%%%%%%%%%%%%%%%%%%%%%%%%%%%%%%%%%%%%%%%%%%%
%%%%%%%%%%%%%%%%%%%%%%%%%%%%%%%%%%%%%%%%%%%%%%%%%%%%%%%%%%%%%%%%%%%%%%%%%%%%%%%%%%%%%%%%%%%%%%%%%%%%

\section{Preliminary results}
\label{sec:survivor}

By a \emph{domain} in \(\CC\) we mean a nonempty, connected set that is open in the spherical topology. 
In particular, \(\CC\) itself is a domain. 
We write \(\sigma\) for the spherical distance on \(\CC\).
For every nonempty set $A\subseteq\mbC$, we write $\diam A:=\sup\{|z-w|:z,w\in A\}$ for its Euclidean diameter.
A family \(\mcF\) of holomorphic maps from a domain \(V\) to \(\CC\) is normal on \(V\) if and only if it is spherically equicontinuous at every point of \(V\) \cite[Theorem~3.2.1]{Schiff1993}.

We first recall Montel's theorem on three fixed omitted values \cite[Theorem~3.7]{Milnor2006}.

\begin{thm}[Montel's theorem]
\label{thm:montel-fixed}
    Let \(V\subseteq\CC\) be a domain and let \(\mcF\) be a family of holomorphic maps \(f\colon V\to\CC\). 
    Suppose there are three distinct points \(a,b,c\in\CC\) such that $f(V)\cap\{a,b,c\}=\emptyset$ for every $f\in\mcF$.
    Then \(\mcF\) is normal on \(V\) with respect to the spherical metric.
\end{thm}

The following extension of Montel's theorem allows the three omitted values to depend on the map, provided they remain uniformly separated in the spherical metric. 
This extension is sometimes attributed to C.~Carathéodory; see \cite{Beardon_2020},
The same statement already appears in \cite[Théorème~J]{Bloch1925}..

We use the following formulation in \cite[Theorem~3.3.5]{Beardon1991}.

\begin{thm}[Generalized Montel's theorem]
\label{thm:montel-separated}
    Let \(V\subseteq\CC\) be a domain and let \(\mcF\) be a family of holomorphic maps \(f\colon V\to\CC\). 
    Suppose there exists \(\delta>0\) such that, for every \(f\in\mcF\), there are three points \(a_f,b_f,c_f\in\CC\setminus f(V)\) satisfying $\min\bigl\{\sigma(a_f,b_f),\,\sigma(b_f,c_f),\,
                   \sigma(c_f,a_f)\bigr\}\geq\delta$.
    Then \(\mcF\) is normal on \(V\) with respect to the spherical metric.
\end{thm}

\medskip 

Let $\seq{f}=(f_n)_{n\geq 1}$ be a sequence of nonconstant rational maps $f_n:\CC\to\CC$.
For each phase $n\geq 1$ and each $k\geq 1$, define
\[
    F_{n,k}\,:=\,f_{n+k-1}\circ\cdots\circ f_n,
\]
and set $F_{n,0}:=\mathrm{id}_{\CC}$.
In this notation, the Fatou and Julia sets at phase \(n\) are
\[
    \fat_n
    = \bigl\{z\in\CC :
       \{F_{n,k}:k\geq1\} \text{ is normal on a neighborhood of } z
      \bigr\},
    \qquad
    \jul_n=\CC\setminus\fat_n.
\]
Here normality is taken with respect to the spherical metric \(\sigma\).
Applying \Cref{thm:montel-separated} to a tail of the forward compositions gives the following criterion.

\begin{cor}
\label{cor:phase-omitted-values}
    Fix \(n\geq 1\) and let \(V\subseteq\CC\) be a domain.
    Suppose there exist \(k_0\geq 1\) and \(\delta>0\) such that, for every \(k\geq k_0\), the image \(F_{n,k}(V)\) omits three points \(a_k,b_k,c_k\in\CC\) satisfying
    \[
        \min\{\sigma(a_k,b_k),\sigma(b_k,c_k),\sigma(c_k,a_k)\}
        \geq \delta.
    \]
    Then \(V\subseteq\fat_n\).
\end{cor}
\begin{proof}
    By \Cref{thm:montel-separated}, the family \(\{F_{n,k}:k\geq k_0\}\) is normal on \(V\). Adding the finitely many earlier compositions preserves normality, so \(V\subseteq\fat_n\).
\end{proof}

We start this section with the following rudimentary result.  

\begin{lem}
\label{lem:pullback-components}
    Let $g:\CC\to\CC$ be a nonconstant rational map and let $W\subseteq \CC$ be a domain i.e. a nonempty, connected open set.
    Then $g^{-1}(W)$ has at most $\deg(g)$ connected components.
    For each such component $U$, the restriction $g|_U:U\to W$ is proper and surjective, and $g(\overline U)=\overline W$.
\end{lem}
\begin{proof}
    Fix a connected component $U$ of $g^{-1}(W)$.
    Since $g^{-1}(W)$ is open in $\CC$, its component $U$ is open in $\CC$ and closed relative to $g^{-1}(W)$.
    For every compact set $K\subseteq W$, the set $g^{-1}(K)$ is compact in $\CC$ and lies in $g^{-1}(W)$.
    Thus $U\cap g^{-1}(K)$ is closed in $g^{-1}(K)$ and hence compact, proving that $g|_U$ is proper.
    
    By the open mapping theorem, $g(U)$ is open in $W$.
    To see that $g(U)$ is also closed in $W$, let $(w_n)_{n\ge 1}$ be a sequence in $g(U)$ that converges to $w\in W$. 
    Then we can find a sequence $(z_n)_{n\ge 1}$ in $U$ so that $g(z_n) = w_n$.
    The set $K := \{w\}\cup\{w_j:j\geq1\}$ is compact in $W$, so properness of $g|_U:U\to W$, established earlier, gives a subsequence, say $(z_{n_k})_{k\ge 1}$ of $(z_n)$ that converges to some $z\in U$.
    By continuity of $g$, we must have that $g(z)=w$.
    Thus, $w\in g(U)$ and therefore $g(U)$ is closed in $W$.
    As $g(U)$ is nonempty and both open and closed in the connected set $W$, we conclude that $g(U)=W$.
    
    By continuity of $g$, we have that $g(\overline U)\subseteq\overline{g(U)}$.
    Conversely, $g(\overline U)$ is compact and hence closed in $\CC$, so it contains the closure of $g(U)$.
    Thus $g(\overline U)=\overline{g(U)}=\overline W$, where the last equality follows from $g(U)=W$.
    
    Finally, fix $w\in W$.
    Since $g(U)=W$ for every component $U$ of $g^{-1}(W)$, each such component contains at least one point of $g^{-1}(\{w\})$.
    The components are pairwise disjoint and $g^{-1}(\{w\})$ contains at most $\deg(g)$ points, so $g^{-1}(W)$ has at most $\deg(g)$ connected components.
\end{proof}

In this section we consider sequences of non-constant rational maps, $\seq{f} = (f_n)_{n\ge 1}$, that satisfy the following condition.  
\begin{equation}
\label{eq:T}
     f_n^{-1}\left( \overline{\mbD}\right)\Subset\mbD\quad(n\geq1),\tag{T}
\end{equation}
We call \eqref{eq:T} the \emph{backward trapping condition}: at every phase, the full preimage of the closed unit disk is compactly contained in the open unit disk. 
Thus, once an orbit leaves \(\mbD\), it cannot return to \(\overline{\mbD}\) at any later phase.
For such sequences of non-constant rational maps, define 
\[
    \mcU_{n,k}
    :=
    \left\{
        \text{connected components of }
        (F_{n,k})^{-1}(\mbD)
    \right\},
    \qquad \text{ for each }\qquad n,k\geq1.
\]
Then using \Cref{lem:pullback-components} we define 
\[
    \Delta_{n,k}:=\max_{U\in\mcU_{n,k}} \diam\overline U.
\]
We impose the following additional assumption on the sequence  $\seq{f} = (f_n)_{n\ge 1}$. 
\begin{equation}
\label{eq:S}
    \Delta_{1,k}\to 0
    \qquad\text{ as }\qquad k\to\infty.\tag{S}
\end{equation}
We call \eqref{eq:S} the \emph{shrinking condition}: as \(k\to\infty\), the largest diameter among the connected components of \(F_{1,k}^{-1}(\mbD)\) tends to zero.

We now introduce the following notation.
For $n\geq 1$ and $k\geq 0$, let
\[
    E_{n,k}:=F_{n,k}^{-1}(\overline{\mbD}),
    \qquad\text{ and }\qquad
    L_n:=\bigcap_{k\geq 0}E_{n,k}.
\]
In particular, $E_{n,0}=\overline{\mbD}$ for every $n\geq 1$, and $L_1=\bigcap_{k\geq 0}E_{1,k}$.

With these notations, we now record a property of the pullback components that we shall use below.
First, recall that a continuous map $g:X\to Y$ is called \emph{proper} if $g^{-1}(A)$ is compact in $X$ for every compact set $A\subseteq Y$.

We now present the central result that shall be useful for the constructions below. 

\begin{thm}
\label{thm:survivor}
    Let $\seq{f} = (f_n)_{n\ge 1}$ be a sequence of non-constant rational maps satisfying \eqref{eq:T}.
    Then the condition \eqref{eq:S} is equivalent to total disconnectedness of \(L_1\). 
    When these equivalent conditions hold, we have that \(\jul_n=L_n\) and \(\Delta_{n,k}\to0\) as $k$ increases, for each phase $n$.
    
    Furthermore, if these conditions hold and \(\#L_1 = \infty\), then every \(\jul_n\) is a Cantor set.
    Under \eqref{eq:S}, every open neighborhood \(V\) of a point of \(\jul_n\) satisfies \(F_{n,k}(V)\supseteq\overline{\mbD}\) for all sufficiently large \(k\).
\end{thm}
\begin{proof}
    First, note that for each fixed $n\in\mbN$, condition \eqref{eq:T} gives
    \begin{equation}
    \label{eq:nesting}
        E_{n,k+1}
        =F_{n,k}^{-1}\bigl(f_{n+k}^{-1}(\overline{\mbD})\bigr)
        \subseteq F_{n,k}^{-1}(\mbD)
        \subseteq E_{n,k}
        \qquad (k\geq 0).       
    \end{equation}
    Thus the sequence $(E_{n,k})_{k\geq 0}$ is nested.
    Since each $F_{n,k}$ is continuous and surjective on $\CC$, every $E_{n,k}$ is nonempty and compact.
    It follows that their intersection $L_n$ is nonempty and compact.
    Moreover, $L_n\subseteq E_{n,1}\subseteq \mbD$.
    
    The sequence $(L_n)_{n\geq 1}$ is invariant under $\seq{f}$ in the sense that
    \begin{equation}
    \label{eq:invariance_of_L_seq}
        L_n=f_n^{-1}(L_{n+1})
        \qquad\text{and}\qquad
        f_n(L_n)=L_{n+1}.
    \end{equation}
    Indeed, the identity $F_{n,k+1}=F_{n+1,k}\circ f_n$ gives $f_n^{-1}(E_{n+1,k})=E_{n,k+1}$ for every $k\geq 0$.
    Consequently,
    \[
        f_n^{-1}(L_{n+1})
        =\bigcap_{k\geq 0}f_n^{-1}(E_{n+1,k})
        =\bigcap_{k\geq 0}E_{n,k+1}
        =L_n,
    \]
    where the last equality follows from the nesting property (\cref{eq:nesting}).
    Applying $f_n$ to both sides and using its surjectivity gives $f_n(L_n)=L_{n+1}$.
    By iteration, we also obtain
    \begin{equation}
    \label{eq:invariance_of_L_seq_2}
        F_{n,k}^{-1}(L_{n+k})=L_n
        \qquad\text{and}\qquad
        F_{n,k}(L_n)=L_{n+k}
        \qquad (k\geq 0).
    \end{equation}
    
    Having observed these properties, we now show that condition \eqref{eq:S} is equivalent to total disconnectedness of $L_1$.
    To this end, assume that \eqref{eq:S} fails.
    Since $\Delta_{1,k}\geq 0$, there exist $\varepsilon>0$ and a strictly increasing sequence $(k_j)_{j\geq 1}$ such that $\Delta_{1,k_j}\geq\varepsilon$ for every $j$.
    By \Cref{lem:pullback-components}, each collection $\mcU_{1,k_j}$ is finite and nonempty, so we may choose $U_j\in\mcU_{1,k_j}$ satisfying $ \diam\overline{U_j}=\Delta_{1,k_j}\geq\varepsilon$.
    Since $\overline{\mbD}$ is compact, a subsequence of $(\overline{U_j})_{j\geq 1}$ converges in the Hausdorff metric to a nonempty compact set $K\subseteq \overline{\mbD}$; \cite[Theorem~3.16]{Falconer85}.
    Each $\overline{U_j}$ is connected, so its Hausdorff limit $K$ is connected as well; \cite[Theorem~3.18]{Falconer85}.
    After relabeling the convergent subsequence, let $h_j := d_H(\overline{U_j},K)$, so that $h_j\to 0$.
    We first show that $\diam K\geq\varepsilon$.
    Since $\overline{U_j}$ is compact, we may choose $x_j,y_j\in\overline{U_j}$ such that $|x_j-y_j|=\diam\overline{U_j}\geq\varepsilon$.
    By the definition of $h_j$ and compactness of $K$, there exist $a_j,b_j\in K$ such that $|x_j-a_j|\leq h_j$ and $|y_j-b_j|\leq h_j$.
    Then, 
    \[
        \varepsilon
        \leq |x_j-y_j|
        \leq |x_j-a_j|+|a_j-b_j|+|b_j-y_j|
        \leq 2h_j+\diam K.
    \]
    Letting $j\to\infty$, we obtain $\diam K\geq\varepsilon$.

    We now show that $K\subseteq L_1$, whence $L_1$ cannot be totally disconnected since $K$ is connected and $\diam(K)\ge \varepsilon>0$. 
    Fix $z\in K$ and choose $z_j\in\overline{U_j}$ with $|z_j-z|\leq h_j$, so that $z_j\to z$.
    For each fixed $k\geq 0$, we have $k_j\geq k$ for all sufficiently large $j$ where we recall that $U_j \in \mcU_{1,k_j}$.
    Since $E_{1,k_j}$ is closed and contains $U_j$, it also contains $\overline{U_j}$.
    Thus, by the nesting property, \ref{eq:nesting},
    \[
        z_j\in\overline{U_j}\subseteq E_{1,k_j}\subseteq E_{1,k}
        \qquad\text{for all sufficiently large }j.
    \]
    Since $E_{1,k}$ is closed, it contains the limit $z$.
    This holds for every $k\geq 0$, so $z\in\bigcap_{k\geq 0}E_{1,k}=L_1$.
    As $z\in K$ was arbitrary, we conclude that $K\subseteq L_1$.

    Conversely, assume that condition \eqref{eq:S} holds, and let $C$ be a nonempty connected subset of $L_1$.
    For every $k\geq 1$, we have $L_1\subseteq E_{1,k+1}\subseteq F_{1,k}^{-1}(\mbD)$, so $C$ lies in a single component of $F_{1,k}^{-1}(\mbD)$.
    Consequently, $\diam C\leq\Delta_{1,k}$ for every $k$.
    Letting $k\to\infty$, we obtain $\diam C=0$, so $C$ is a singleton.
    This proves that $L_1$ is totally disconnected whenever \eqref{eq:S} holds

    \medskip

    Now for the rest of the proof, assume that \eqref{eq:S} holds or equivalently that $L_1$ is totally disconnected. 

    Fix a phase $n\in\mbN$ and put $H := F_{1,n-1}$, where $H$ is the indenity map if $n=1$.
    For each $k\geq 1$ and each $U\in\mcU_{n,k}$, choose a connected component $V$ of $H^{-1}(U)$.
    By \Cref{lem:pullback-components}, we have $H(V)=U$ and $H(\overline V)=\overline U$.
    We claim that $V$ is also a connected component of $F_{1,n+k-1}^{-1}(\mbD)$.
    Indeed, since $F_{1,n+k-1}=F_{n,k}\circ H$, the set $V$ lies in some component $V'$ of $F_{1,n+k-1}^{-1}(\mbD)$.
    Then $H(V')$ is a connected subset of $F_{n,k}^{-1}(\mbD)$ meeting $U$, so $H(V')\subseteq U$.
    Thus $V\subseteq V'\subseteq H^{-1}(U)$, and the maximality of $V$ gives $V'=V$.
    
    Since $V$ is a component of $F_{1,n+k-1}^{-1}(\mbD)$, we have $\diam\overline V\leq\Delta_{1,n+k-1}$.
    To see that $\Delta_{n,k}\to 0$, fix $\varepsilon>0$.
    The set $E_{1,n-1}=H^{-1}(\overline{\mbD})$ is a compact subset of $\overline{\mbD}$, and $H$ is continuous on this set with values in $\overline{\mbD}$.
    Thus $H$ is uniformly continuous on $E_{1,n-1}$ in the Euclidean metric.
    Choose $\delta>0$ such that $|H(z)-H(w)|<\varepsilon$ whenever $z,w\in E_{1,n-1}$ and $|z-w|<\delta$.
    By \eqref{eq:S}, there exists $k_0$ such that $\Delta_{1,n+k-1}<\delta$ for every $k\geq k_0$.
    For each such $k$ and each $U\in\mcU_{n,k}$, the corresponding $\overline V$ lies in $E_{1,n-1}$ because $H(\overline V)=\overline U\subseteq\overline{\mbD}$.
    Since $\diam\overline V<\delta$, it follows that $\diam\overline U=\diam H(\overline V)\leq\varepsilon$.
    Taking the maximum over $U\in\mcU_{n,k}$ gives $0\le \Delta_{n,k}\leq\varepsilon$ for every $k\geq k_0$, proving that $\Delta_{n,k}\to 0$.
    
    Now let $C$ be a nonempty connected subset of $L_n$.
    For every $k\geq 1$, we have $C\subseteq L_n\subseteq E_{n,k+1}\subseteq F_{n,k}^{-1}(\mbD)$, so $C$ lies in a single component of $F_{n,k}^{-1}(\mbD)$.
    Consequently, $\diam C\leq\Delta_{n,k}$ for every $k$.
    Letting $k\to\infty$, we obtain $\diam C=0$, so $C$ is a singleton.
    Hence $L_n$ is totally disconnected.

    \medskip
    
    We next show that $\jul_n=L_n$.
    If $z\notin L_n$, then $z\notin E_{n,k_0}$ for some $k_0\geq 0$.
    Since $E_{n,k_0}$ is closed, there is an open neighborhood $V$ of $z$ disjoint from $E_{n,k_0}$.
    By the nesting property (\ref{eq:nesting}), $V$ is disjoint from $E_{n,k}$ for every $k\geq k_0$, so each $F_{n,k}$ with $k\geq k_0$ omits $\overline{\mbD}$ on $V$.
    Montel's theorem, \Cref{thm:montel-fixed} (and \Cref{cor:phase-omitted-values}), now gives normality of the family on $V$, and hence $z\in\fat_n$.
    This proves $\jul_n\subseteq L_n$.
    
    Now let $z\in L_n$ and let $V$ be any open neighborhood of $z$.
    For each $k\geq 1$, the inclusion $L_n\subseteq F_{n,k}^{-1}(\mbD)$ gives a unique component $U_k\in\mcU_{n,k}$ containing $z$.
    Since $\diam\overline{U_k}\leq\Delta_{n,k}\to 0$ as $k$ increases, we have $\overline{U_k}\subseteq V$ for all sufficiently large $k$.
    By \Cref{lem:pullback-components}, for every such sufficiently large $k$ we have
    \begin{equation}
    \label{eq:covering_property}
         F_{n,k}(V)\supseteq F_{n,k}(\overline{U_k})=\overline{\mbD}.
    \end{equation}
    Suppose, for a contradiction, that \(z\notin\jul_n\). 
    Then the family \(\{F_{n,k}:k\geq0\}\) is normal on a neighborhood of \(z\) and hence spherically equicontinuous at \(z\). 
    Taking \(\varepsilon := \sigma(0,1)/3\), we can choose an open neighborhood \(V\)
    of \(z\) such that
    \[
        \sigma\bigl(F_{n,k}(w),F_{n,k}(z)\bigr)<\varepsilon
        \qquad(w\in V,\ k\geq0).
    \]
    By \eqref{eq:covering_property} above, \(F_{n,k}(V)\) contains both \(0\) and \(1\) for all sufficiently large \(k\). For any such \(k\),we have that 
    \[
        \sigma(0,1)
        \leq \sigma\bigl(0,F_{n,k}(z)\bigr)+\sigma\bigl(F_{n,k}(z),1\bigr)
        <2\varepsilon=\tfrac23\sigma(0,1),
    \]
    a contradiction. 
    Thus \(L_n\subseteq\jul_n\) concluding that \(\jul_n=L_n\) simultaneously with the and the covering property \eqref{eq:covering_property} at every point of \(\jul_n\).

    \medskip 
    
    Finally, suppose that $L_1$ is infinite.
    We have already established that $\jul_n = L_n$, for each $n\in\mbN$ and that $L_n$ is non-empty and compact.
    We also saw that condition \eqref{eq:S} implies that $L_n$ is totally disconnected for each $n\in\mbN$.
    All that remains to show is that $L_n$ and therefore $\jul_n$ is perfect. 
    Since $L_1=F_{1,m-1}^{-1}(L_m)$ by \eqref{eq:invariance_of_L_seq} and as $F_{1,m-1}$ has finite fibers, $L_m$ is infinite for every $m\geq 1$.
    Fix $n\in\mbN$ and fix $z\in L_n$ and let $U_k$ be the component containing $z$ as above.
    For each $k\geq 1$, choose $w_k\in L_{n+k}$ distinct from $F_{n,k}(z)$, this is possible as each $L_m$ is infinite.
    Since $w_k\in L_{n+k}\subseteq E_{n+k,1} = f_{n+k}^{-1}(\overline{\mbD})$, we have from the condtion \eqref{eq:T} that $w_k \in \mbD$.
    Furthermore, from \Cref{lem:pullback-components} we also have that $F_{n,k}(U_k)=\mbD$, thus there exists $z_k\in U_k$ such that $F_{n,k}(z_k)=w_k$.
    But \eqref{eq:invariance_of_L_seq_2} gives $F_{n,k}^{-1}(L_{n+k})=L_n$ and in turn we get that $z_k\in L_n$, and the choice of $w_k$ gives $z_k\neq z$.
    Since $z$ and $z_k$ both belong to $U_k$, we have $|z_k-z|\leq\Delta_{n,k}\to 0$, as $k\to\infty$.
    Therfore $z_k\to z$ and $z_k\in L_n \setminus \{z\}$ for every $k$ whence we have that every neighborhood of \(z\) meets \(L_n\setminus\{z\}\).
    Therefore $L_n = \jul_n$ is also perfect. 
    Combining all properties, we see that each iterated Julia set $\jul_n$ is a Cantor set. 
\end{proof}

\begin{rem}
\label{rem:area_positive_for_all}
    For any sequence of nonconstant rational maps, positivity of planar area at one phase implies positivity at every phase.
    Indeed, applying \eqref{thm:invariance deterministic} to a tail gives \(\jul_n=f_n^{-1}(\jul_{n+1})\) and \(f_n(\jul_n)=\jul_{n+1}\).
    A nonconstant rational map and its local inverse branches preserve null sets off their finitely many exceptional points, so \(\area(\jul_n)=0\) if and only if \(\area(\jul_{n+1})=0\).
\end{rem}

We have the following immediate corollary.

\begin{cor}
\label{cor:area}
    Let $\seq{f}$ be a sequence of nonconstant rational maps satisfying \eqref{eq:T} and \eqref{eq:S}.
    Then
    \[
        \area(\jul_1)
        =\lim_{k\to\infty}\area(E_{1,k})
        =\inf_{k\geq 0}\area(E_{1,k}).
    \]
    If this infimum is positive, then $\jul_n$ is a Cantor set of positive planar area for every $n\in\mbN$.
\end{cor}
\begin{proof}
    By \Cref{thm:survivor}, we have $\jul_1=L_1=\bigcap_{k\geq 0}E_{1,k}$.
    Since the sets $E_{1,k}$ are nested and $\area(E_{1,0})=\area(\overline{\mbD})=\pi<\infty$, the area identity follows from continuity of measure from above.
    Suppose now that $\inf_{k\geq 0}\area(E_{1,k})>0$.
    Then $L_1=\jul_1$ has positive planar area and is therefore infinite.
    It follows from \Cref{thm:survivor} that $\jul_n$ is a Cantor set for every $n\in\mbN$.
    Since $\area(\jul_1)>0$, \Cref{rem:area_positive_for_all} gives $\area(\jul_n)>0$ for every $n\in\mbN$.
\end{proof}

By \Cref{cor:area}, any lower bound $\area(E_{1,k})\geq a>0$ valid for every $k\geq 1$ gives $\area(\jul_1)\geq a$.
Such a bound can also be obtained by estimating the area lost at each step.
Indeed, since $E_{1,k}\subseteq E_{1,k-1}$ and $\area(E_{1,0})=\pi$, summing the successive losses gives
\[
    \area(\jul_1)
    =\pi-\sum_{k=1}^{\infty}\area(E_{1,k-1}\setminus E_{1,k}).
\]
Thus, if $\area(E_{1,k-1}\setminus E_{1,k})\leq\delta_k$ for every $k\geq 1$ and $\sum_{k=1}^{\infty}\delta_k<\pi$, then $\area(\jul_1)\geq\pi-\sum_{k=1}^{\infty}\delta_k>0$.

\begin{rem}
\label{rem:trapping-disk}
    The unit disk in \Cref{thm:survivor,cor:area} may be replaced throughout by any fixed open disk $D := B(c,R)$ with $R>0$.
    Indeed, let $A(z) := (z-c)/R$ and put $\widetilde f_n := A\circ f_n\circ A^{-1}$.
    Then $f_n^{-1}(\overline D)\subseteq D$ if and only if $\widetilde f_n^{-1}(\overline{\mbD})\Subset\mbD$.
    Since $\widetilde F_{n,k}=A\circ F_{n,k}\circ A^{-1}$, the map $A$ takes each component of $F_{n,k}^{-1}(D)$ onto a component of $\widetilde F_{n,k}^{-1}(\mbD)$ and divides its Euclidean diameter by $R$.
    Thus the shrinking condition is equivalent in the two coordinates.
    The phase Julia sets correspond under $A$, and the covering conclusion of \Cref{thm:survivor} becomes $F_{n,k}(V)\supseteq\overline D$ for all sufficiently large $k$.
\end{rem}

%%%%%%%%%%%%%%%%%%%%%%%%%%%%%%%%%%%%%%%%%%%%%%%%%%%%%%%%%%%%%%%%%%%%%%%%%%%%%%%%%%%%%%%%%%%%%%%%%%%%
%%%%%%%%%%%%%%%%%%%%%%%%%%%%%%%%%%%%%%%%%%%%%%%%%%%%%%%%%%%%%%%%%%%%%%%%%%%%%%%%%%%%%%%%%%%%%%%%%%%%

\section{Proof of Theorem~A}
\label{sec:constructions}

We apply \Cref{thm:survivor,cor:area} to prove \Cref{thm:main}.
For each sequence under consideration, recall that $F_{n,k}$ denotes the composition of $k$ successive maps
starting at phase $n$, with $F_{n,0}=\mathrm{id}_{\CC}$. We write
\[
    E_{n,k} := F_{n,k}^{-1}(\overline{\mbD}),
    \qquad
    L_n := \bigcap_{k\geq 0}E_{n,k}.
\]
Thus $L_n$ consists of the points whose forward orbits, starting at phase $n$, stay in $\overline{\mbD}$. 
In particular, for the polynomial sequence $(f_n)_{n\geq 1}$ constructed below, we have
\[
    E_{n,1}=f_n^{-1}(\overline{\mbD}),
    \qquad
    E_{1,n}=(f_n\circ\cdots\circ f_1)^{-1}(\overline{\mbD}).
\]

%%%%%%%%%%%%%%%%%%%%%%%%%%%%%%%%%%%%%%%%%%%%%%%%%%%%%%%%%%%%%%%%%%%%%%%%%%%%%%%%%%%%%%%%%%%%%%%%%%%%
%%%%%%%%%%%%%%%%%%%%%%%%%%%%%%%%%%%%%%%%%%%%%%%%%%%%%%%%%%%%%%%%%%%%%%%%%%%%%%%%%%%%%%%%%%%%%%%%%%%%

\subsection{Existence of the polynomial sequence}
\label{sec:polynomial}

We prove the existence of the polynomial sequence in part~(1) of \Cref{thm:main} inductively.

\begin{thm}
\label{thm:polynomial}
    For every \(a\in(0,\pi)\), there exists a sequence \(\seq{f}=(f_n)_{n\geq1}\) of polynomials in \(\mbQ[z]\), each of degree at least two, such that \(\seq{f}\) satisfies \eqref{eq:T} and, for every \(n\geq1\),
    \[
        \area(E_{n,1})>1,
        \qquad
        a<\area(E_{1,n})<a+2^{-n}(\pi-a),
        \qquad\text{ and }\qquad 
        \Delta_{1,n}<2^{-n}.
    \]
    In particular, \(\area(\jul_1)=a\), and for every \(n\geq1\), the set \(\jul_n\) is a Cantor set of positive planar area compactly contained in \(\mbD\).
    Moreover, \(\area(\jul_n\cap B(z,r))>0\) for every \(n\geq1\), every \(z\in\jul_n\), and every \(r>0\).
\end{thm}
\begin{proof}
    Fix $0<a<\pi$.
    We construct the polynomials $f_n$ inductively, starting with $F_{1,0}=\mathrm{id}_{\CC}$ and $E_{1,0}=\overline{\mbD}$, so that $\area(E_{1,0})=\pi>a$.
    At step $n$, suppose that $f_1,\ldots,f_{n-1}$ have been chosen and that $\area(E_{1,n-1})>a$.
    Put $F := F_{1,n-1}$, so that $F=\mathrm{id}_{\CC}$ when $n=1$.
    Write $F(z)=cz^D+\cdots + c_0$, where $D := \deg F$ and $c,c_0\in\mbQ$ with $c\neq 0$.
 
    We first estimate the diameters of the connected components of $F^{-1}(B)$ for a nonempty set $B\subseteq \mbC$ with $\diam B\leq\delta$.
    We claim that every such component $C$ satisfies
    \begin{equation}
    \label{eq:poly-inverse-bound}
        \diam C\leq 2D\left(\frac{\delta}{|c|}\right)^{1/D}.
    \end{equation}
    To see this, fix $w\in B$ and write $F(z)-w=c\prod_{j=1}^D(z-\xi_j)$.
    If $z\in F^{-1}(B)$, then $|F(z)-w|\leq\delta$, so $|z-\xi_j|\leq r:=(\delta/|c|)^{1/D}$ for at least one $j\in \{1, \ldots, D\}$.
    Thus $F^{-1}(B)$ is contained in the union of the $D$-many closed disks of radius $r$ centered at the points $\xi_j$.
    Any two points in the same connected component of this union lie in disks joined by a chain of at most $D$ disks, with consecutive disks intersecting.
    Their distance is therefore at most $2Dr$, which proves \eqref{eq:poly-inverse-bound}.
    Now let
    \[
        \delta := 
        \max\{|\operatorname{Re}c|,|\operatorname{Im}c|\} \times 
        \left(\frac{2^{-n}}{4D}\right)^D.
    \]
    Since $\max\{|\operatorname{Re}c|,|\operatorname{Im}c|\}\leq |c|$, the right-hand side of \eqref{eq:poly-inverse-bound} is at most $2^{-n}/2$.

    Now we define a finite Borel measure, $\mu$, on $\overline{\mbD}$ as follows. 
    For any Borel set $B\subseteq\overline{\mbD}$, we define $\mu(B) := \area(F^{-1}(B))$.
    This defines a finite measure because $F^{-1}(\overline{\mbD})=E_{1,n-1}$ is a compact subset of $\mbC$.
    Moreover, $\mu$ is absolutely continuous with respect to planar area, since the preimage of a set of area zero under a nonconstant polynomial has area zero.
    In particular, $\mu(\partial\mbD)=0$, and hence
    \[
        \mu(\mbD)=\mu(\overline{\mbD})=\area(E_{1,n-1})>a.
    \]
    Since the imaginary axis divides $\mbD$ into two half-disks. 
    But as \( \mu(\mbD)\) is positive, we must have that one of the half-disks must have $\mu$-measure at most \(\mu(\mbD)/2\). 

    Suppose first that the left half-disk has \(\mu\)-measure at most \(\mu(\mbD)/2\).
    For each \(t\in[0,1]\), put
    \[
        \Omega_t := \{z\in\mbD:\operatorname{Re}z<t\}.
    \]
    Then \(\Omega_0\) is the left half-disk and \(\Omega_1=\mbD\), so
    \[
        \mu(\Omega_0)\leq\frac{\mu(\mbD)}{2}
        \qquad\text{and}\qquad
        \mu(\Omega_1)=\mu(\mbD).
    \]
    The function \(t\mapsto\mu(\Omega_t)\) is continuous because every vertical line has \(\mu\)-measure zero.
    Since \(0<a<\mu(\mbD)\), the interval
    \[
        \left(
            \max\left\{a,\frac{\mu(\mbD)}{2}\right\},
            \frac{a+\mu(\mbD)}{2}
        \right)
    \]
    is non-empty and lies between $\mu(\Omega_0)$ and $\mu(\Omega_1)$.
    We may therefore choose \(t_0\in(0,1)\cap\mbQ\) such that
    \begin{equation}
    \label{ineq:measure-of_omega-t}
        \max\left\{a,\frac{\mu(\mbD)}{2}\right\}
        <\mu(\Omega_{t_0})
        <\frac{a+\mu(\mbD)}{2}.     
    \end{equation}
    If instead the right half-disk has \(\mu\)-measure at most \(\mu(\mbD)/2\), put
    \[
        \Omega_t := \{z\in\mbD:\operatorname{Re}z>-t\}
        \qquad(t\in[0,1]).
    \]
    Now \(\Omega_0\) is the right half-disk and \(\Omega_1=\mbD\), so the same argument gives a rational \(t_0\in(0,1)\) satisfying \eqref{ineq:measure-of_omega-t}.
    In either case, set \(\Omega := \Omega_{t_0}\).
    The set \(\Omega\) contains the chosen half-disk, so \(\area(\Omega)\geq\pi/2>1\) and also symmetric about the real axis. 
    
    For each $m\ge 1$, let $A_m$ be the union of all closed dyadic squares of side length $2^{-m}$ contained in $\Omega$.
    Since $A_m\uparrow\Omega$, continuity from below for measures, gives $\mu(A_m)\to\mu(\Omega)>a$ and $\area(A_m)\to\area(\Omega)>1$. Choose $m$ sufficiently large so that
    \[
        \sqrt{2}\,2^{-m}<\delta,\qquad
        \mu(A_m)>a,\qquad \text{ and }\qquad 
        \area(A_m)>1.
    \]
    Denote the squares at this level by $Q_1,\ldots,Q_N$.

    We now use the collection $Q_1,\ldots,Q_N$ to obtain smaller open squares $U_j$ whose closures are pairwise disjoint. 
    Inside each $U_j$ we keep a closed square $K_j$, chosen so that their union still has $\mu$-measure greater than $a$ and planar area greater than $1$.

    Let $s := 2^{-m}$ be the common side length of the squares $Q_j$.
    For each rational $\rho\in(0,1)$, let $K_j(\rho)$ be the closed square concentric with $Q_j$, with sides parallel to those of $Q_j$ and side length $\rho s$. 
    As $\rho\uparrow 1$, the squares $K_j(\rho)$ increase to $\operatorname{int}(Q_j)$. 
    The boundary of each $Q_j$  has planar area zero and hence also $\mu$-measure zero by absolute continuity. 
    Thus, writing
    \[
        K(\rho) := \bigcup_{j=1}^N K_j(\rho)
        \qquad\text{and}\qquad Q := \bigcup_{j=1}^N Q_j,
    \]
    we obtain, by continuity from below of measures, that
    \[
        \mu(K(\rho)) \to \mu(Q)>a
        \qquad\text{and}\qquad
        \area(K(\rho)) \to \area(Q)>1
        \qquad\text{ as } \quad \rho\uparrow 1.
    \]
    Choose a rational $\rho$ sufficiently close to $1$ that both strict inequalities hold for $K := K(\rho)$, and set $K_j := K_j(\rho)$.
    For each $j$, let $U_j$ be the concentric open square with sides parallel to those of $Q_j$ and side length $(1+\rho)s/2$. 
    Since $\rho s < (1+\rho)s/2 < s$ we have 
    \[
        K_j\subseteq U_j\subseteq \overline{U_j}
        \subseteq \operatorname{int}(Q_j).
    \]
    The center of $Q_j$ has rational coordinates and the side length of $U_j$ is rational, so the vertices of $U_j$ lie in $\mbQ(i)$.
    Since the interiors of the $Q_j$ are pairwise disjoint, the closures of the $U_j$ are pairwise disjoint as well. 
    Finally, $\diam U_j<\diam Q_j<\delta$ for every $j$, and the finite union $U := \bigcup_{j=1}^N U_j$ is compactly contained in $\mbD$.
    
    Since $K$ is a finite union of pairwise disjoint closed squares, its complement in $\mbC$ is connected.
    Since $U$ is an open neighborhood of $K$, Hilbert's lemniscate theorem \cite[Theorem~5.5.8]{Ransford_1995} gives a polynomial $\tilde{q}$ such that
    \[
        \dfrac{|\tilde{q}(z)|}{\|\tilde{q}\|_K} > 1 \qquad\text{ for }\qquad z\in \mbC\setminus U.
    \]
    Here, $\|\tilde{q}\|_K := \max_{w\in K}|\tilde{q}(w)|$.
    Now take $q :=  \tilde{q} / \|\tilde{q}\|_K$ then $|q|\leq 1$ on $K$ and $|q|>1$ on $\mbC\setminus U$.
    The polynomial $q$ is nonconstant, so $|q(z)|\to\infty$ as $|z|\to\infty$.
    It follows that $|q|$ attains its minimum $b$ on the closed set $\mbC\setminus U$ and that $b>1$.
    Since $b>1$, take $t := (1/b+1)/2$, so that $1/b<t<1$.
    Then
    \[
        \max_{w\in K}|tq(w)|\leq t<1<tb
        =\min_{w\in\mbC\setminus U}|tq(w)|.
    \]

    Having obtained a polynomial $q$ with complex coefficients, we now approximate the coefficients of $tq$ by Gaussian rationals while preserving the two strict inequalities above. 
    Write
    \[
        tq(z)=\sum_{j=0}^d a_j z^j,\qquad \text{ with }\qquad a_d\neq 0.
    \]
    We first choose a radius $R$ beyond which the leading term dominates for all sufficiently small perturbations of the coefficients. 
    Indeed, if $p(z) := \sum_{j=0}^d b_jz^j$, with $|b_d-a_d|<|a_d|/2$ and $|b_j-a_j|<1$ for $j<d$, then for $|z|\geq 1$ we have
    \[
        |p(z)|\geq \frac{|a_d|}{2}|z|^d
            -\left(\sum_{j=0}^{d-1}(|a_j|+1)\right)|z|^{d-1}.
    \]
    Put $M := \sum_{j=0}^{d-1}(|a_j|+1)$ and choose $R>1$ so large that
    \[
        \frac{|a_d|}{2}R-M>1.
    \]
    This choice of $R$ depends only on the coefficients of $tq$. 
    For every polynomial $p$ satisfying the above coefficient bounds and every $z$ with $|z|\geq R$, the preceding estimate gives
    \[
        |p(z)|
        \geq |z|^{d-1}\left(\frac{|a_d|}{2}|z|-M\right)
        \geq \frac{|a_d|}{2}R-M
        >1,
    \]
    where we used $|z|^{d-1}\geq 1$.
    
    It remains to control $p$ on $|z|\leq R$. Set
    \[
        \eta := \frac12\min\{1-t,tb-1\}>0.
    \]
    By density of $\mbQ(i)$ in $\mbC$, choose $b_0,\ldots,b_d\in\mbQ(i)$ sufficiently close to $a_0 \ldots,a_d$ that the preceding coefficient bounds hold and
    \[
        \sum_{j=0}^d |b_j-a_j|R^j<\eta.
    \]
    Then $|p(z)-tq(z)|<\eta$ for every $|z|\leq R$. 
    Since $K\subseteq \mbD$, we obtain $|p(z)|\leq t+\eta<1$ on $K$.
    For $z\in\mbC\setminus U$ with $|z|\leq R$, we likewise obtain $|p(z)|\geq tb-\eta>1$. Together with the bound for $|z|\geq R$, this gives $|p|>1$ on all of $\mbC\setminus U$.
    Finally, $|b_d-a_d|<|a_d|/2$ ensures that $b_d\neq 0$, so $p\in\mbQ(i)[z]$ has the same degree as $q$.

    The dyadic grid and \(\Omega\) are symmetric about the real axis, so reflection across this axis permutes the squares forming \(A_m\).
    The concentric squares used to form \(K\) and \(U\) therefore have this symmetry, as well. 
    Put \(p^*(z) := \overline{p(\bar z)}\) and set \(f_n := p\,p^*\).
    Since \(p\in\mbQ(i)[z]\), the coefficients of \(f_n\) belong to \(\mbQ(i)\).
    Conjugating the coefficients interchanges \(p\) and \(p^*\), leaving their product unchanged.
    Thus each coefficient of \(f_n\) is rational, and \(f_n\in\mbQ[z]\).
    Moreover, \(|p^*(z)|=|p(\bar z)|\), so the symmetry of \(K\) and \(U\) gives \(|f_n|<1\) on \(K\) and \(|f_n|>1\) outside \(U\).
    So,
    \[
         K\subseteq f_n^{-1}(\mbD)
        \subseteq E_{n,1}=f_n^{-1}(\overline{\mbD})
        \subseteq U\Subset\mbD.
    \]
    Each $U_j$ meets $f_n^{-1}(\mbD)$ because it contains $K_j$. Since $f_n^{-1}(\mbD)\subseteq U$ and the $U_j$ are pairwise disjoint open squares, each component of $f_n^{-1}(\mbD)$ lies in a single $U_j$. 
    Thus $f_n^{-1}(\mbD)$ has at least $N$ components.
    By \Cref{lem:pullback-components}, $\deg(f_n)\geq N$. In particular, $\deg(f_n)\geq 2$, since $N4^{-m}=\area(A_m)>1$.
    Thus $f_n$ satisfies \eqref{eq:T} and $\area(E_{n,1})\geq\area(K)>1$.
    Since $F_{1,n}=f_n\circ F$, we have $E_{1,n}=F^{-1}(E_{n,1})$ and hence
    \[
        a<\mu(K)\leq\area(E_{1,n})=\mu(E_{n,1})\leq\mu(\Omega)<(a+\mu(\mbD))/2.
    \]

    To estimate $\Delta_{1,n}$, fix $V\in\mcU_{1,n}$.
    Since $V\subseteq E_{1,n}$ and $E_{1,n}$ is closed, we have $\overline V\subseteq E_{1,n}$.
    Thus $F(\overline V)$ is a connected subset of $E_{n,1}\subseteq U$ and lies in a single square $U_j$.
    It follows that $\overline V$ lies in a single connected component of $F^{-1}(U_j)$.
    Since $\diam U_j<\delta$, estimate \eqref{eq:poly-inverse-bound} gives $\diam\overline V\leq 2^{-n-1}$.
    Taking the maximum over $V\in\mcU_{1,n}$ yields $\Delta_{1,n}\leq 2^{-n-1}<2^{-n}$.

    This constructs $f_n$ with all the required properties, and the inequality $\area(E_{1,n})>a$ allows us to proceed to step $n+1$.
    By induction, we obtain the required sequence $\seq{f}$.
    Since $\Delta_{1,n}<2^{-n}$ for every $n$, condition \eqref{eq:S} holds.
    Since $\mu(\mbD) =\area(E_{1,n-1})$, iteration gives $0<\area(E_{1,n})-a<2^{-n}(\pi-a)$.
    By \Cref{cor:area}, we have $\area(\jul_1)=\lim_{n\to\infty}\area(E_{1,n}) = a$, and every $\jul_n$ is a Cantor set of positive planar area.
    Then, \Cref{thm:survivor} gives $\jul_n=L_n\subseteq E_{n,1}\Subset\mbD$ for every $n\geq 1$.

    Ruffianly,fix $n\geq 1$, $z\in\jul_n$, and $r>0$. 
    By \Cref{thm:survivor}, we may choose $k$ such that $F_{n,k}(B(z,r))\supseteq\overline{\mbD}$. 
    Since $\jul_{n+k}\subseteq\mbD$ and $F_{n,k}^{-1}(\jul_{n+k})=\jul_n$, it follows that $F_{n,k}\bigl(B(z,r)\cap\jul_n\bigr)=\jul_{n+k}$. 
    Hence $\area(B(z,r)\cap\jul_n)>0$, as all iterated Julia sets $\jul_m$ must have positive area from \Cref{rem:area_positive_for_all} and as polynomials map planar null sets to planar null sets. 
\end{proof}

\begin{cor}
\label{cor:polynomail}
     For every $a>0$, there exists a sequence of polynomials $f_n\in\mbQ[z]$ with $\deg(f_n)\geq 2$ for every $n\geq 1$ such that $\area(\jul_1)=a$.
     Furthermore, there is an open disk $D$ centered at the origin such that, for every $n\geq 1$, we have $f_n^{-1}(\overline D)\Subset D$, and $\jul_n$ is a Cantor set compactly contained in $D$. 
     Moreover,
    \[
        \area\bigl(\jul_n\cap B(z,r)\bigr)>0
        \qquad\text{for every }z\in\jul_n\text{ and every }r>0.
    \]
\end{cor}
\begin{proof}
    Fix \(a>0\) and choose a positive integer \(L\) such that \(b:=a/L^2<\pi\).
    Put \(S(z):=Lz\) and \(D:=S(\mbD)=B(0,L)\).
    By \Cref{thm:polynomial}, there exists a sequence of polynomials \(f_n\in\mbQ[z]\), each of degree at least two, satisfying \eqref{eq:T} and such that \(\area(\jul_1)=b\).
    Each \(\jul_n\) is a Cantor set compactly contained in \(\mbD\) and has positive planar area in every neighborhood of each of its points.

    Define \(\widetilde f_n := S\circ f_n\circ S^{-1}\).
    Since \(L\) is a positive integer, \(\widetilde f_n(z)=L f_n(z/L)\) belongs to \(\mbQ[z]\) and has the same degree as \(f_n\).
    The forward compositions satisfy \(\widetilde F_{n,k}=S\circ F_{n,k}\circ S^{-1}\).
    Since fixed conjugation preserves normality, the phase Julia sets satisfy \(\widetilde{\jul}_n=S(\jul_n)\).
    Thus each \(\widetilde{\jul}_n\) is a Cantor set compactly contained in \(D\), and
    \[
        \widetilde f_n^{-1}(\overline D)
        =S\bigl(f_n^{-1}(\overline{\mbD})\bigr)\Subset D.
    \]
    Dilation by \(L\) multiplies planar area by \(L^2\), so \(\area(\widetilde{\jul}_1)=L^2b=a\).
    Moreover, for every \(z\in\jul_n\) and \(r>0\),
    \[
        \area\bigl(\widetilde{\jul}_n\cap B(S(z),Lr)\bigr)
        =L^2\area\bigl(\jul_n\cap B(z,r)\bigr)>0.
    \]
    This proves the polynomial assertion.
\end{proof}
\begin{rem}
    A \emph{bounded polynomial sequence} in the sense of Comerford~\cite[\S1.1]{Comerford2013} has degrees between $2$ and a fixed upper bound. Its coefficients are uniformly bounded in modulus, and its leading coefficients are bounded away from zero.
    The polynomial sequence constructed above is \emph{not} a sequence of bounded polynomials since $\deg(f_n)\to\infty$.
    Indeed, from the proof above we have that \(\operatorname{diam} Q_j<\delta\), the covering used in \eqref{eq:poly-inverse-bound} shows that \(F^{-1}(Q_j)\) lies in the union of \(D\) disks of radius at most \(2^{-n}/(4D)\). 
    Summing their areas over the \(N\) squares gives
    \[
        a<\mu(A_m)\le ND\pi\left(\frac{2^{-n}}{4D}\right)^2
        =\frac{\pi N}{16D}\,4^{-n}.
    \]
    Hence \(\deg(f_n)\ge N>(16a/\pi)D\,4^n\). Since \(D\ge1\), the degrees tend to infinity.
\end{rem}

%%%%%%%%%%%%%%%%%%%%%%%%%%%%%%%%%%%%%%%%%%%%%%%%%%%%%%%%%%%%%%%%%%%%%%%%%%%%%%%%%%%%%%%%%%%%%%%%%%%%
%%%%%%%%%%%%%%%%%%%%%%%%%%%%%%%%%%%%%%%%%%%%%%%%%%%%%%%%%%%%%%%%%%%%%%%%%%%%%%%%%%%%%%%%%%%%%%%%%%%

\subsection{Rational map construction}
\label{sec:rational}

Conjugating the polynomial sequence of \Cref{thm:polynomial} by a suitable fixed M\"obius automorphism of $\mbD$ already gives a sequence of rational maps with finite poles whose phase Julia sets are Cantor sets of positive planar area. 
Part~(2) of \Cref{thm:main} also rules out polynomial sequences obtained by changing M\"obius coordinates separately at each phase.
The following theorem specifies such a sequence recursively. We prescribe the area in \Cref{cor:prescribed-area}.

\begin{thm}
\label{thm:rational}
    There exists a sequence $\seq{f}=(f_n)_{n\geq1}$ of rational maps in $\mbQ(z)$ satisfying \eqref{eq:T} such that no $f_n$ is M\"obius-conjugate to a polynomial and, for every $n\geq1$,
    \[
        \Delta_{1,n}<2^{-n-1}
        \qquad\text{and}\qquad
        \area(E_{1,n-1}\setminus E_{1,n})<2^{-n-1}.
    \]
    In particular, \(\area(\jul_1)>0\).
    There is no sequence $(M_n)_{n\geq1}$ of M\"obius maps for which $M_{n+1}\circ f_n\circ M_n^{-1}$ is a polynomial for every $n\geq1$.
    For every $n\geq1$, the phase Julia set $\jul_n$ is a Cantor set with $\jul_n\Subset\mbD$, and $\area(\jul_n\cap B(z,r))>0$ for every $z\in\jul_n$ and every $r>0$.
\end{thm}

\begin{proof}
     We begin with $F_{1,0}=\mathrm{id}$ and $E_{1,0}=\overline{\mbD}$.
    At step $n$, suppose that $f_1,\ldots,f_{n-1}$ have been chosen and satisfy \eqref{eq:T}.
    Put $F := F_{1,n-1}$ and write $F=A/B$, where $A$ and $B$ are relatively prime polynomials and $B$ is monic.
    At the first step, $A(z)=z$ and $B(z)=1$.
    In every case, $F^{-1}(\overline{\mbD})=E_{1,n-1}\subseteq\overline{\mbD}$, by \eqref{eq:T} at the preceding steps.
    In particular, $F(\infty)\notin\overline{\mbD}$.

    Put $D := \deg F$ and write $A(z)=\sum_{j=0}^D a_jz^j$ and $B(z)=\sum_{j=0}^D b_jz^j$, allowing zero coefficients.
    Define
    \[
        b_*:=\sum_{j=0}^D|b_j|,
        \qquad\text{and}\qquad
        c:=|a_D|-|b_D|.
    \]
    We have \(b_*>0\) because \(B\) is monic.
    If \(b_D\neq0\), then \(F(\infty)=a_D/b_D\notin\overline{\mbD}\), so \(|a_D|>|b_D|\).
    If \(b_D=0\), then \(\deg B<D\), and hence \(\deg A=D\), so again \(|a_D|>|b_D|\).
    Thus \(c>0\).

    Now let \(H:=A'B-AB'\).
    Since $F$ is nonconstant, $H$ is a nonzero polynomial.
    Let $m := \deg H$ and let $\gamma$ be its leading coefficient.
    Set
    \[
        e:=2^{-n},\qquad
        r:=\left(\frac{e}{4\pi(m+1)}\right)^{1/2},\qquad
        \kappa:=\frac{|\gamma|r^m}{b_*^2},\qquad
        \tau:=\min\left\{\frac1{100},\frac{e\kappa^2}{4D}\right\}.
    \]
    Take \(\beta:=1/\lceil32/\tau\rceil\).
    Thus \(\beta\) is rational and \(0<\beta\leq\tau/32\).
    Take
    \begin{equation}
    \label{eq:r-grid}
        N:=1+\left\lceil
            \max\left\{4,\frac{16\pi}{\tau},
            \frac{2b_*}{c}(2^{n+2}D)^D\right\}
        \right\rceil. 
    \end{equation}

    Take as centers the points in
    \[
        \Lambda_N:=\{(j+ik)/N:j,k\in\mbZ,\ j^2+k^2\leq(N-1)^2\}.
    \]
    For each \(u\in\Lambda_N\), take three squares centered at \(u\) with sides parallel to the coordinate axes.
    Let \(K_u\) be the closed inner square, \(\Gamma_u\) the boundary of the middle square, and \(U_u\) the open outer square.
    Their respective side lengths are
    \[
        \frac{1-\beta}{N},\qquad
        \frac{1-\beta/2}{N},\qquad
        \frac{1-\beta/4}{N}.
    \]
    Put \(K := \bigcup_{u\in\Lambda_N}K_u\) and \(U := \bigcup_{u\in\Lambda_N}U_u\).
    Distinct centers differ by at least \(1/N\) in at least one coordinate, so the outer squares have pairwise disjoint closures.
    Moreover, \(U\Subset\mbD\), since \(|u|\leq1-1/N\) and each outer square has half-diagonal less than \(1/N\).

    Now we are ready to define the next function $f_n$. 
    Put
    \[
        M:=\left\lceil\frac{1024N^2}{\beta^2}\right\rceil.
    \]
    For each \(u\in\Lambda_N\), orient the boundary \(\Gamma_u\) of the middle square counterclockwise and divide each of its four sides into \(M\) equal segments.
    Starting at the lower left vertex, number these segments from \(1\) to \(4M\) in the direction of the orientation.
    Let \(S_{u,\nu}\) denote the oriented segment from its initial endpoint \(a_{u,\nu}\) to its terminal endpoint \(b_{u,\nu}\), and put
    \[
        \zeta_{u,\nu}:=\frac{a_{u,\nu}+b_{u,\nu}}{2},
        \qquad
        d_{u,\nu}:=b_{u,\nu}-a_{u,\nu}.
    \]
    Since \(\Lambda_N\) is finite, we may define \(f_n\) by the finite sum
    \begin{equation}
    \label{eq:r-map}
        f_n(z):=2-\frac{1}{3i}
        \sum_{u\in\Lambda_N}\sum_{\nu=1}^{4M}
        \frac{d_{u,\nu}}{\zeta_{u,\nu}-z}.
    \end{equation}
    
    The centers and vertices have rational real and imaginary parts, and the same holds for the endpoints and midpoints of the segments.
    Thus \(\zeta_{u,\nu},d_{u,\nu}\in\mbQ(i)\), so \(f_n\in\mbQ(i)(z)\).
    Reflection in the real axis reverses orientation and therefore pairs each midpoint and increment \((\zeta,d)\) with \((\bar\zeta,-\bar d)\) in the sum.
    Since \(\overline{1/(3i)}=-1/(3i)\), it follows that \(f_n(\bar z)=\overline{f_n(z)}\).
    By uniqueness of the reduced quotient with monic denominator, its coefficients are fixed by conjugation and hence lie in \(\mbQ\).
    Thus \(f_n\in\mbQ(z)\).
    We also record that rotation by a quarter-turn permutes the oriented segments and sends \((\zeta,d)\) to \((i\zeta,id)\).
    Reindexing the finite number of summand in \eqref{eq:r-map} therefore gives \(f_n(iz)=f_n(z)\).
    The midpoints are distinct, and every increment is nonzero, so each \(\zeta_{u,\nu}\) is a simple pole of \(f_n\).
    The formula also gives \(f_n(\infty)=2\).

    Having defined \(f_n\), we now verify that 
    \[
        \Delta_{1,n}<2^{-n-1}
        \qquad\text{and}\qquad
        \area(E_{1,n-1}\setminus E_{1,n})<2^{-n-1}.
    \]
    We first show that
    \[
        K\subseteq f_n^{-1}(\mbD)
        \subseteq E_{n,1}=f_n^{-1}(\overline{\mbD})
        \subseteq U.
    \]

    For \(z\) off the middle-square boundaries, let
    \[
        I(z):=\frac{1}{2\pi i}
        \sum_{u\in\Lambda_N}\int_{\Gamma_u}\frac{d\zeta}{\zeta-z},
        \qquad\text{ and }\qquad
        R(z):=\frac{1}{2\pi i}
        \sum_{u\in\Lambda_N}\sum_{\nu=1}^{4M}
        \frac{d_{u,\nu}}{\zeta_{u,\nu}-z}.
    \]
    If \(z\in K\), then \(z\) lies in one inner square \(K_u\), hence inside the corresponding middle square and outside all the others.
    By Cauchy's integral formula and Cauchy's theorem, \(I(z)=1\).
    If \(z\in\mbC\setminus U\), then \(z\) lies outside every middle square, so \(I(z)=0\).
    
    Fix \(z\in K\cup(\mbC\setminus U)\).
    The gap between corresponding sides of an inner square and its middle square is \(\beta/(4N)\), and the gap between corresponding sides of a middle square and its outer square is \(\beta/(8N)\).
    For each \(u\in\Lambda_N\), either \(z\in K_u\) or \(z\notin U_u\), since the outer squares are pairwise disjoint.
    It follows that the distance from \(z\) to every \(\Gamma_u\) is at least \(\beta/(8N)\).
    
    Fix \(u\in\Lambda_N\) and \(1\leq\nu\leq4M\).
    Recall that the segment \(S_{u,\nu}\) denote the oriented segment from its initial endpoint \(a_{u,\nu}\) to its terminal endpoint \(b_{u,\nu}\), and has midpoint \(\zeta_{u,\nu}\) and length \(L := (1-\beta/2)/(NM)\).
    For every \(\zeta\in S_{u,\nu}\), we have
    \[
        \left|\frac{1}{\zeta-z}-\frac{1}{\zeta_{u,\nu}-z}\right|
        =
        \frac{|\zeta-\zeta_{u,\nu}|}
             {|\zeta-z|\,|\zeta_{u,\nu}-z|}
        \leq
        \left(\frac{8N}{\beta}\right)^2|\zeta-\zeta_{u,\nu}|.
    \]
    Integrating along \(S_{u,\nu}\) gives
    \[
    \begin{aligned}
        \left|\int_{S_{u,\nu}}\frac{d\zeta}{\zeta-z}
              -\frac{d_{u,\nu}}{\zeta_{u,\nu}-z}\right|
        &=
        \left|\int_{S_{u,\nu}}
            \left(\frac{1}{\zeta-z}-\frac{1}{\zeta_{u,\nu}-z}\right)d\zeta\right|\\
        &\leq
        \left(\frac{8N}{\beta}\right)^2
        \int_{S_{u,\nu}}|\zeta-\zeta_{u,\nu}|\,|d\zeta|\\
        &=
        \left(\frac{8N}{\beta}\right)^2
        2\int_0^{L/2}t\,dt
        =
        \left(\frac{8N}{\beta}\right)^2\frac{L^2}{4}.
    \end{aligned}
    \]
    There are at most \((2N-1)^2<4N^2\) centers and \(4M\) segments per center.
    Summing the preceding estimate over all segments, dividing by \(2\pi\), and using \(L<1/(NM)\), we obtain
    \[
        |R(z)-I(z)|
        \leq \frac{16MN^2}{2\pi}
              \left(\frac{8N}{\beta}\right)^2\frac{L^2}{4}
        \leq \frac{128N^2}{\pi M\beta^2}
        \leq \frac{1}{8\pi}.
    \]

    By \eqref{eq:r-map}, \(f_n(z)=2-(2\pi/3)R(z)\).
    For \(z\in K\), we have \(I(z)=1\), so
    \[
        |f_n(z)|
        \leq 2(\pi/3-1)+\frac{2\pi}{3}|R(z)-1|
        \leq 2(\pi/3-1)+\frac{1}{12}
        <\frac14.
    \]
    For \(z\in\mbC\setminus U\), we have \(I(z)=0\), so
    \[
        |f_n(z)-2|=\frac{2\pi}{3}|R(z)|\leq\frac{1}{12}.
    \]
    Thus \(|f_n|<1\) on \(K\) and \(|f_n|>1\) on \(\mbC\setminus U\).
    Together with \(f_n(\infty)=2\), this proves the claimed inclusions.
    Since \(U\Subset\mbD\), condition \eqref{eq:T} holds at step \(n\).

    Now we proceed to show that $\Delta_{1,n}<2^{-n-1}$.
    Fix \(u\in\Lambda_N\) and consider its outer square \(U_u\).
    If \(z\in F^{-1}(U_u)\), then $z\in E_{1,n-1}\subseteq\overline{\mbD}$, so $|B(z)|\leq b_*$ and
    \[
        |A(z)-uB(z)|=|B(z)|\,|F(z)-u|\leq\frac{2b_*}{N}.
    \]
    The polynomial $A-uB$ has degree $D$ and its leading coefficient has modulus at least $c$, since $|u|<1$.
    The argument in \eqref{eq:poly-inverse-bound} therefore gives, for every connected component \(C\) of \(F^{-1}(U_u)\),
    \[
        \diam C
        \leq 2D\left(\frac{2b_*}{cN}\right)^{1/D}
        <2^{-n-1},
    \]
    where the strict inequality follows from \eqref{eq:r-grid}.
    Now let $V\in\mcU_{1,n}$.
    Since $\overline V$ is connected and lies in $E_{1,n}=F^{-1}(E_{n,1})$, its image under $F$ is a connected subset of $E_{n,1}\subset U$.
    Hence $F(\overline V)$ lies in a single outer square \(U_u\), and $\overline V$ lies in a single connected component of \(F^{-1}(U_u)\).
    The preceding bound is independent of $V$, so $\Delta_{1,n}<2^{-n-1}$.

    It remains to estimate $\area(E_{1,n-1}\setminus E_{1,n})$.
    The closed squares of side \(1/N\) about the chosen centers cover \(\{|z|\leq1-2/N\}\), since a nearest grid point to any such \(z\) lies at distance at most \(1/(\sqrt{2}N)\) and has modulus at most \(1-1/N\).
    Replacing each such square by \(K_u\) removes at most \(2\beta/N^2\) of its area.
    Since there are fewer than $4N^2$ squares, we obtain
    \[
        \area(\overline{\mbD}\setminus K)
        \leq\pi\bigl(1-(1 - 2/N)^2\bigr)+8\beta
        <4\pi/N+8\beta
        <\frac{\tau}{2}.
    \]

    Let $P$ be the union of the open disks of radius $r$ centered at the zeros of $H$, with $P=\emptyset$ if $H$ is constant.
    Then $\area(P)\leq m\pi r^2<e/4$.
    For $z\in\overline{\mbD}\setminus P$ with $B(z)\neq 0$, factoring $H$ at its zeros gives
    \[
        |F'(z)|=\frac{|H(z)|}{|B(z)|^2}
        \geq\frac{|{\gamma}|r^m}{b_*^2}
        =\kappa.
    \]
    Put $S := F^{-1}(\overline{\mbD}\setminus K)$.
    Since $K\subseteq E_{n,1}$, we have $E_{1,n-1}\setminus E_{1,n}\subseteq S\subseteq E_{1,n-1}$.
    The set $S$ contains no pole of $F$, and each value of $F$ has at most $D$ preimages.
    Applying the area formula to $F$ on $S\setminus P$ yields
    \[
        \kappa^2\area(S\setminus P)
        \leq\int_{S\setminus P}|F'(z)|^2\,d\area(z)
        \leq D\,\area(\overline{\mbD}\setminus K).
    \]
    It follows from the choice of $\tau$ that
    \[
        \area(E_{1,n-1}\setminus E_{1,n})
        \leq\area(P)+\area(S\setminus P)
        <\frac e4+\frac{D\tau}{2\kappa^2}
        \leq\frac{3e}{8}
        <2^{-n-1}.
    \]

    This constructs $f_n$ with all the required properties, and we proceed to the next step with $F_{1,n}=f_n\circ F$.
    By induction, we obtain the sequence $\seq{f}$, and the bound on $\Delta_{1,n}$ gives \eqref{eq:S}.
    Since $E_{1,0}=\overline{\mbD}$ and the sets $E_{1,n}$ are nested, \Cref{cor:area} gives
    \[
        \area(\jul_1)
        =\pi-\sum_{n=1}^{\infty}\area(E_{1,n-1}\setminus E_{1,n})
        \geq\pi-\sum_{n=1}^{\infty}2^{-n-1}
        =\pi-\frac{1}{2}.
    \]
    The same corollary shows that every $\jul_n$ is a Cantor set of positive planar area.
    Finally, $\jul_n=L_n\subseteq E_{n,1}\Subset\mbD$ by \Cref{thm:survivor} and \eqref{eq:T}.
    
    Fix \(n\geq1\).
    We first show that \(0\) has at least two distinct preimages under \(f_n\).
    Each inner square \(K_u\) lies in \(f_n^{-1}(\mbD)\), and every component of \(f_n^{-1}(\mbD)\) lies in a single outer square, since \(f_n^{-1}(\mbD)\subseteq U\) and the outer squares have pairwise disjoint closures.
    Thus inner squares with different centers lie in different components of \(f_n^{-1}(\mbD)\).
    There are at least two such components, and by \Cref{lem:pullback-components} each maps onto \(\mbD\) under \(f_n\).
    Each therefore contains a preimage of \(0\).

    We now show that \(f_n\) is not M\"obius-conjugate to a polynomial.
    Suppose that \(f_n=T\circ P\circ T^{-1}\) for some polynomial \(P\) and some M\"obius map \(T\).
    Since \(P^{-1}(\{\infty\})=\{\infty\}\), the point \(p := T(\infty)\) satisfies
    \[
        f_n^{-1}(\{p\})=\{p\}.
    \]
    In particular, \(f_n(p)=p\).
    Recall that \(f_n(iz)=f_n(z)\) by the rotational symmetry of the grid.
    If \(p\) were finite and nonzero, then \(ip\neq p\) and \(f_n(ip)=f_n(p)=p\), giving a second preimage of \(p\).
    Also, \(p\neq\infty\) because \(f_n(\infty)=2\), and \(p\neq0\) because \(0\) has at least two preimages.
    This proves the claim.

    Suppose now that there was a sequence of M\"obius maps \((M_n)_{n\ge 1}\) such that \(M_{n+1}\circ f_n\circ M_n^{-1}\) were a polynomial for every \(n\geq1\).
    Putting \(x_n := M_n^{-1}(\infty)\), we would have
    \[
        f_n^{-1}(\{x_{n+1}\})=\{x_n\}
        \qquad(n\geq1).
    \]
    If \(x_n\) were finite and nonzero, then \(ix_n\) would be a second preimage of \(x_{n+1}\).
    Hence \(x_n\in\{0,\infty\}\) at every phase.
    If \(x_n=\infty\) for some \(n\), then \(x_{n+1}=f_n(\infty)=2\), contrary to \(x_{n+1}\in\{0,\infty\}\).
    It follows that \(x_n=0\) for every \(n\), which would give \(f_n^{-1}(\{0\})=\{0\}\), again a contradiction.
    
    Finally, fix \(n\geq1\), \(z\in\jul_n\), and \(r>0\).
    By \Cref{thm:survivor}, \(F_{n,k}(B(z,r))\) contains \(\overline{\mbD}\) for all sufficiently large \(k\).
    Since \(\jul_{n+k}\subseteq\mbD\) and \(F_{n,k}^{-1}(\jul_{n+k})=\jul_n\), it follows that
    \[
        F_{n,k}\bigl(B(z,r)\cap\jul_n\bigr)=\jul_{n+k}.
    \]
    The image \(\jul_{n+k}\) has positive planar area, and \(F_{n,k}\) is locally Lipschitz away from its poles, all of which lie outside \(\jul_n\).
    Hence \(\area(B(z,r)\cap\jul_n)>0\).
\end{proof}

We now prescribe the planar area of the first phase Julia set.

\begin{thm}
\label{cor:prescribed-area}
    For every $a>0$, there exists a sequence of rational maps $f_n\in\mbQ(z)$, none of which is M\"obius-conjugate to a polynomial, such that $\area(\jul_1)=a$.
    This rational sequence admits no phase-dependent M\"obius change of coordinates to a polynomial sequence, in the sense of \Cref{thm:rational}.
    Furthermore, there is an open disk $D$ centered at the origin such that, for every $n\geq 1$, we have $f_n^{-1}(\overline D)\Subset D$, and $\jul_n$ is a Cantor set compactly contained in $D$. 
    Moreover,
    \[
        \area\bigl(\jul_n\cap B(z,r)\bigr)>0
        \qquad\text{for every }z\in\jul_n\text{ and every }r>0.
    \]
\end{thm}
\begin{proof}
    Fix \(a>0\) and choose a positive integer \(L\) such that \(b:=a/L^2<\pi\).
    Put \(S(z) := Lz\) and \(D := S(\mbD)=B(0,L)\).

    We first construct a sequence \((g_n)_{n\geq1}\) in \(\mbQ(z)\) whose first phase Julia set has planar area \(b\).
    At step \(n\), we use the proof of \Cref{thm:rational} to construct a grid map \(f_n\), and then choose a rational number \(t_n>1\) and set \(g_n := t_nf_n\).
    
    Put \(G_0 := \id\) and \(E_0 := \overline{\mbD}\).
    After choosing \(g_1,\ldots,g_{n-1}\), write
    \[
        G_{n-1} := g_{n-1}\circ\cdots\circ g_1,
        \qquad
        E_{n-1} := G_{n-1}^{-1}(\overline{\mbD}).
    \]
    Suppose that the maps chosen so far satisfy \eqref{eq:T} and that \(\area(E_{n-1})>b\).
    At the first step this area inequality holds because \(\area(E_0)=\pi>b\).
    
    Choose an integer \(j\geq n\) large enough that
    \[
        2^{-j-1}<\area(E_{n-1})-b.
    \]
    Apply the grid construction from the proof of \Cref{thm:rational} to \(G_{n-1}\), using \(j\) as the index in its parameters, so that \(e=2^{-j}\).
    The construction applies because \(E_{n-1}\subseteq\overline{\mbD}\) and \(G_{n-1}(\infty)\notin\overline{\mbD}\).
    Denote the resulting grid map by \(f_n\in\mbQ(z)\) and put \(R_n := f_n\circ G_{n-1}\).
    By the estimates in that proof, \(f_n^{-1}(\overline{\mbD})\Subset\mbD\) and
    \[
        \area\bigl(E_{n-1}\setminus R_n^{-1}(\overline{\mbD})\bigr)
        <2^{-j-1},
        \qquad
        \diam\overline V<2^{-j-1}
        \quad\text{for every component \(V\) of \(R_n^{-1}(\mbD)\).}
    \]
    Since \(R_n^{-1}(\overline{\mbD})\subseteq E_{n-1}\), the first estimate and the choice of \(j\) give
    \[
        \area\bigl(R_n^{-1}(\overline{\mbD})\bigr)
        >\area(E_{n-1})-2^{-j-1}>b.
    \]
    
    We now choose \(t_n\) so that the surviving area is still greater than \(b\), but exceeds \(b\) by less than \(2^{-n}\).
    If we take \(g_n=tf_n\), then the new composition is \(tR_n\), and a point \(z\) survives this step precisely when \(|R_n(z)|\leq1/t\).
    Accordingly, for \(t\geq1\), put
    \[
        H_n(t) := \area\bigl(\{z:|R_n(z)|\leq1/t\}\bigr).
    \]
    We have just shown that \(H_n(1)>b\).
    The sublevel sets lie in \(E_{n-1}\), which has finite area, and each level set \(\{|R_n|=1/t\}\) has area zero because it is the preimage of a circle under a nonconstant rational map.
    Thus \(H_n\) is continuous by dominated convergence.
    As \(t\to\infty\), the sublevel sets decrease to the finite zero set of \(R_n\), so \(H_n(t)\to0\).
    
    Choose a real number \(y\) satisfying
    \[
        b<y<\min\{H_n(1),\,b+2^{-n}\}.
    \]
    By the intermediate value theorem, \(H_n(t)=y\) for some \(t>1\).
    By continuity, \(b<H_n(t')<b+2^{-n}\) for every \(t'\) sufficiently close to \(t\).
    We may therefore choose a rational \(t_n>1\) satisfying these strict inequalities.
    Set
    \[
        g_n := t_nf_n,\qquad
        G_n := g_n\circ G_{n-1}=t_nR_n,\qquad
        E_n := G_n^{-1}(\overline{\mbD}).
    \]
    Then \(g_n\in\mbQ(z)\), and the choice of \(t_n\) gives
    \[
        b<\area(E_n)=H_n(t_n)<b+2^{-n}.
    \]
    In particular, \(\area(E_n)>b\), as required for the next step.
    
    We also retain the trapping and shrinking estimates.
    Indeed, since \(t_n>1\),
    \[
        g_n^{-1}(\overline{\mbD})
        \subseteq f_n^{-1}(\overline{\mbD})\Subset\mbD.
    \]
    Moreover, the closure of each component of \(G_n^{-1}(\mbD)\) is a connected subset of
    \[
        E_n=\{|R_n|\leq1/t_n\}\subseteq R_n^{-1}(\mbD).
    \]
    It therefore lies in a single component of \(R_n^{-1}(\mbD)\) and has diameter less than \(2^{-j-1}\leq2^{-n-1}\).
    This completes the induction.
    
    Write \(\jul_n\) for the phase-\(n\) Julia set of the sequence \((g_n)_{n\geq1}\).
    This sequence satisfies \eqref{eq:T} and \eqref{eq:S}, and \(E_n\) is its \(n\)-step pullback of \(\overline{\mbD}\) starting at phase one.
    By \Cref{cor:area},
    \[
        \area(\jul_1)=\lim_{n\to\infty}\area(E_n)=b,
    \]
    where the last equality follows from \(0<\area(E_n)-b<2^{-n}\).
    The same corollary and \Cref{thm:survivor} show that each \(\jul_n\) is a Cantor set of positive planar area compactly contained in \(\mbD\).
    
    For the local area assertion, fix an open disk \(B\) meeting \(\jul_n\), and write \(G_{n,k} := g_{n+k-1}\circ\cdots\circ g_n\).
    By the covering property in \Cref{thm:survivor} and phase invariance,
    \[
        G_{n,k}(B\cap\jul_n)=\jul_{n+k}
        \qquad\text{for all sufficiently large \(k\).}
    \]
    Since the image has positive planar area and a rational map takes planar null sets to planar null sets, \(\area(B\cap\jul_n)>0\).
    
    At each step, the grid map \(f_n\) satisfies \(f_n(iz)=f_n(z)\), has at least two distinct zeros, and satisfies \(f_n(\infty)=2\).
    Multiplication by \(t_n\) preserves the rotation identity and the zeros, and gives \(g_n(\infty)=2t_n\), which is finite and nonzero.
    The two conjugacy assertions now follow by the same arguments as in the proof of \Cref{thm:rational}.
    
    Finally, apply the fixed dilation \(S(z)=Lz\) from the polynomial case.
    The maps \(\widetilde g_n := S\circ g_n\circ S^{-1}\) belong to \(\mbQ(z)\), since \(\widetilde g_n(z)=L g_n(z/L)\) and \(L\) is a positive integer.
    Their phase Julia sets are \(S(\jul_n)\), their trapping disk is \(B(0,L)\), and their first phase Julia set has planar area \(L^2b=a\).
    Dilation also preserves the local positive-area property and both conjugacy obstructions.
\end{proof}

%%%%%%%%%%%%%%%%%%%%%%%%%%%%%%%%%%%%%%%%%%%%%%%%%%%%%%%%%%%%%%%%%%%%%%%%%%%%%%%%%%%%%%%%%%%%%%%%%%%%
%%%%%%%%%%%%%%%%%%%%%%%%%%%%%%%%%%%%%%%%%%%%%%%%%%%%%%%%%%%%%%%%%%%%%%%%%%%%%%%%%%%%%%%%%%%%%%%%%%%%

\subsection{Positive-area Cantor Julia sets without a common escape condition}
\label{sec:exterior}

We say that a sequence of polynomials \((f_n)\) satisfies the  \emph{common exterior escape condition} if there exist \(R>0\) and \(\lambda>1\) such that
\begin{equation}
\label{eq:uniform-exterior}
    |f_n(z)|\geq\lambda |z| \qquad(n\geq1,\ z\in\mbC,\ |z|\geq R).
\end{equation}
For polynomials of degree at least two, this condition is equivalent to the hypothesis on a common invariant domain used by B\"uger~\cite{selfsimilar}.

\begin{prop}
\label{rem:buger-common-escape}
    Let \((f_n)\) be a sequence of polynomials with \(\deg f_n\geq2\) for every \(n\geq1\). 
    Then \eqref{eq:uniform-exterior} holds if and only if there is a hyperbolic domain \(M\subseteq\CC\) containing \(\infty\) such that \(f_n(M)\subseteq M\) for every \(n\geq1\) and \(F_{1,k}\to\infty\) locally uniformly on \(M\) as \(k\to\infty\).
\end{prop}
\begin{proof}
    Suppose \eqref{eq:uniform-exterior} holds and take \(M := \{z\in\mbC:|z|>R\}\cup\{\infty\}\).
    This is a hyperbolic domain and each \(f_n\) maps it into itself.
    Iterating the inequality gives \(|F_{1,k}(z)|\geq\lambda^k|z|\) for every finite \(z\in M\). Since each polynomial fixes \(\infty\), it follows that \(F_{1,k}\to\infty\) locally uniformly on \(M\).

    Conversely, suppose such a domain \(M\) exists. 
    Its complement \(L := \CC\setminus M\) is a compact subset of \(\mbC\) and contains two distinct points \(a,b\) by hyperbolicity. 
    Choose \(r>0\) with \(L\subseteq\{|z|\leq r\}\) and set \(\delta := |b-a|\).
    Since \(f_n^{-1}(L)\subseteq L\), all the roots \(\zeta_{n,j}\) of \(f_n-a\) and every preimage of \(b\) lie in this disk.
    Write \(d_n := \deg f_n\) and let \(c_n\) be its leading coefficient.
    Evaluating \(f_n-a\) at any \(w_n\in f_n^{-1}(\{b\})\) gives
    \[
        \delta=|c_n|\prod_{j=1}^{d_n}|w_n-\zeta_{n,j}|
        \leq |c_n|(2r)^{d_n}.
    \]
    Thus for \(|z|\geq3r\),
    \[
        |f_n(z)|\geq |c_n|(|z|-r)^{d_n}-|a|
        \geq \delta\left(\frac{|z|-r}{2r}\right)^2-|a|,
    \]
    where the last inequality uses \(d_n\geq2\) and \((|z|-r)/(2r)\geq1\).
    The last expression exceeds \(2|z|\) for all sufficiently large \(|z|\), uniformly in \(n\). This proves \eqref{eq:uniform-exterior}.
\end{proof}

The polynomial sequence \(\seq{f} = (f_n)_{n\ge 1}\) of \Cref{thm:polynomial} satisfies the two equivalent conditions in \Cref{rem:buger-common-escape}. 
Indeed, by \eqref{eq:T}, the function \(w\mapsto1/f_n(1/w)\) extends across \(0\) to a holomorphic self-map of \(\mbD\) with a zero of order \(d_n=\deg f_n\) at \(0\).
Schwarz's lemma gives 
\begin{equation}
\label{eq:lemniscate-exterior}
    |f_n(z)|\geq |z|^{d_n} \qquad (|z|>1).
\end{equation}
Since \(d_n\geq2\), we can take \(R=\lambda=2\) in \eqref{eq:uniform-exterior}.

The next proposition gives a polynomial sequence for which neither of these equivalent conditions holds.  
For every \(R>0\), some rescaled polynomial sends a point of \(\{|z|>R\}\) into \(\{|z|<R\}\). The rescaling preserves the first
Julia set exactly, and every iterated Julia set is still a compact planar Cantor set of positive area.

\begin{prop}
\label{prop:exterior}
    Fix \(0<a<\pi\), and let \(\seq{f}=(f_n)_{n\geq1}\) be a  polynomial sequence from \Cref{thm:polynomial}, with \(d_n := \deg f_n\). Denote its phase-\(n\) Julia set by \(\jul_n\).
    Put \(A_0:=1\) and, for \(n\geq1\), define
    \[
        A_n:=
        \begin{cases}
            2^{n d_n}, & n \text{ even},\\
            1,         & n \text{ odd}.
        \end{cases}
    \]
    Define the rescaled sequence \(\seq{\widetilde f}\) by
    \[
        \widetilde f_n(z) := A_n^{-1}f_n(A_{n-1}z),
    \]
    and denote its phase-\(n\) Julia set by \(\widetilde\jul_n\).

    Then, for every \(R>0\), some \(\widetilde f_n\) sends a point of \(\{|z|>R\}\) into \(\{|z|<R\}\). At every phase \(n\), the set \(\widetilde\jul_n\) is a compact planar Cantor set of positive area, and
    \[
        \widetilde\jul_n=A_{n-1}^{-1}\jul_n,
        \qquad
        \area(\widetilde\jul_n)=A_{n-1}^{-2}\area(\jul_n).
    \]
    In particular, \(\widetilde\jul_1=\jul_1\) and \(\area(\widetilde\jul_1) = a\). 
    Moreover, for every \(n\geq1\), \(z\in\widetilde\jul_n\), and \(r>0\),
    \[
        \area\bigl(\widetilde\jul_n\cap B(z,r)\bigr)>0.
    \]
\end{prop}
\begin{proof}
    Write \(f_n(z)=c_n\prod_{j=1}^{d_n}(z-\zeta_{n,j})\). 
    Recall that \(E_{n,1}=f_n^{-1}(\overline{\mbD})\). 
    By \Cref{thm:polynomial}, \(E_{n,1}\Subset\mbD\) and \(\area(E_{n,1})>1\). 
    Each root \(\zeta_{n,j}\) belongs to \(E_{n,1}\). 
    Therefore, 
    \begin{align*}
        \int_{E_{n,1}}\max\{0,-\log|z-\zeta_{n,j}|\}\,d\area(z)
            &
                =\int_0^{2\pi}\int_0^1
                    \mathbf 1_{E_{n,1}}(\zeta_{n,j}+te^{i\theta})\,(-\log t)\,t\,dt\,d\theta\\
            &\leq 2\pi\int_0^1(-\log t)\,t\,dt
            =\frac{\pi}{2}.
    \end{align*}

    Now, using the factorization of $f_n$ above, we have that 
    \[
        \log{|f_n(z)|} = \log{|c_n|} + \sum_{j=1}^{d_n} \log{|z - \zeta_{n,j}|}
    \]
    On $E_{n,1}$, the left-hand side above is at most $0$.
    Therefore, on $E_{n,1}$ we have that 
    \[
        \log{|c_n|}
        \le \sum_{j=1}^{d_n} - \log{|z - \zeta_{n,j}|}.
    \]
    Integrating over \(E_{n,1}\), we then get that 
    \[
        \area{(E_{n,1})} \log{|c_n|}
            \le \sum_{j=1}^{d_n} \int_{E_{n,1}}\max\{0, -\log{|z - \zeta_{n,j}|}\} \, d\area(z)
             \le \sum_{j=1}^{d_n}  \dfrac{\pi}{2} = \dfrac{\pi d_n}{2}. 
    \]
    But from \Cref{thm:polynomial} we know that $\area(E_{n,1})> 1$, therefore $|c_n|\leq e^{\pi d_n/2}$. 

    Now let $n$ be an even positive integer.
    Thus, \(A_{n-1}=1\) and \(A_n=2^{n d_n}\). 
    Therefore, 
    \[
        \widetilde f_n(z) = A_n^{-1} f_n(A_{n-1}z) = 2^{-nd_n}f_n(z). 
    \]
    Now for fixrd $r>0$ and for any $|z|\le r$ we have that $|z- \zeta_{n,j}| \le r+ 1$ and therefore, for the even $n$ we have 
    \[
        |\widetilde f_n(z)| \le 2^{-nd_n} e^{\pi d_n/ 2} (r+1)^{d_n}. 
    \]
    This means that     
    \[
        \sup_{|z|\leq r}|\widetilde f_n(z)|
        \xrightarrow[\substack{n\to\infty\\ n\text{ even}}]{} 0.
    \]

    Now to prove the first assertion, let $R>0$.
    Take $r := R+1$ then we have that  \(|\widetilde f_n(R+1)|<R\) for all sufficiently large \emph{even} \(n\),  proving the first assertion with $z = R+1 \in \{|w| > R\}$.

    Next we prove that  $ \widetilde\jul_n=A_{n-1}^{-1}\jul_n$, from which the assertion that $\area(\widetilde\jul_n)=A_{n-1}^{-2}\area(\jul_n)$, follows immediately. 

    Fix a phase \(n\), and put \(\mu_n(z) := A_{n-1}z\). For \(k\geq1\),  write
    \[
        F_{n,k}=f_{n+k-1}\circ\cdots\circ f_n,
        \qquad\text{ and similary }\qquad 
        \widetilde F_{n,k}
         := \widetilde f_{n+k-1}\circ\cdots\circ\widetilde f_n,
    \]
    with \(\widetilde F_{n,0}:=\mathrm{id}=F_{n,0}\). 
    Then we have 
    \[
        \widetilde F_{n,k}=H_{n,k}\circ\mu_n,
        \qquad\text{ where }\quad 
        H_{n,k}(u) := A_{n+k-1}^{-1}F_{n,k}(u).
    \]
    We first identify the Julia set of the family \((H_{n,k})_{k\geq1}\).
    By \Cref{thm:survivor}, applied to the original sequence $\seq{f}$, we get that
    \[
        \jul_n=\bigcap_{k\geq1}E_{n,k},
        \qquad\text{ where }\qquad 
        E_{n,k}=F_{n,k}^{-1}(\overline{\mbD}).
    \]
    Hence for every \(u\in\jul_n\) and \(k\geq1\),
    \[
        |H_{n,k}(u)|
        =A_{n+k-1}^{-1}|F_{n,k}(u)|
        \leq A_{n+k-1}^{-1}\leq1.
    \]

    Now fix \(u\in\mbC\setminus\jul_n\). 
    There there is some \(q\geq1\) for which \(u\notin E_{n,q}\), so \(|F_{n,q}(u)|>1\).
    By continuity, there is a neighborhood \(U\) of \(u\) and there also exists a number \(\rho>1\) such that for all $v\in U$ we have that $|F_{n,q}(v)|\geq\rho$. 
    By \eqref{eq:lemniscate-exterior}, each subsequent step satisfies
    \[
        |F_{n,s+1}(v)|
        =|f_{n+s}(F_{n,s}(v))|
        \geq |F_{n,s}(v)|^2
        \qquad(s\geq q).
    \]
    Iterating this inequality \(k-1-q\) times gives $ |F_{n,k-1}(v)|\geq\rho^{\,2^{k-1-q}}$, for all $v\in U$ and for each $k\ge q+1$. 

    Put \(m := n+k-1\), so that \(F_{n,k}=f_m\circ F_{n,k-1}\).
    By \eqref{eq:lemniscate-exterior} and the bound \(A_m\leq 2^{m d_m}\), we have, for every \(v\in U\) and \(k\geq q+1\),
    \[
        |H_{n,k}(v)|
        =\frac{|f_m(F_{n,k-1}(v))|}{A_m}
        \geq \frac{|F_{n,k-1}(v)|^{d_m}}{2^{m d_m}}
        \geq
        \left(\frac{\rho^{\,2^{k-1-q}}}{2^{\,n+k-1}}\right)^{d_m}.
    \]
    Since \(\rho>1\) and \(n,q\) are fixed, and as $\rho >1$ and \(d_m\geq2\), it follows that \(H_{n,k}\to\infty\) uniformly on \(U\).
    Here notice that the numerator \(\rho^{2^{k-1-q}}\) grows doubly exponentially in \(k\), whereas the denominator \(2^{n+k-1}\) grows exponentially. 
    For all finite \(v\) with \(|v|\geq2\), the same bound holds with \(q=0\) and \(\rho=2\), since \(F_{n,0}(v)=v\).
    Together with \(H_{n,k}(\infty)=\infty\), this gives locally uniform spherical convergence to \(\infty\) on a neighborhood of \(\infty\).
    Thus \((H_{n,k})_{k\geq1}\) is normal on a neighborhood of every \(z\in\CC\setminus\jul_n\). Hence $\jul\left((H_{n,k})_{k\geq1}\right)\subseteq\jul_n$.

    For the reverse inclusion, let \(B\) be an open disk meeting \(\jul_n\).
    Since \(\jul_n\) is a compact Cantor set, \(B\setminus\jul_n\) is a nonempty open subset of \(B\). 
    If the family \((H_{n,k})_{k\geq1}\) were normal on \(B\), there would exist a subsequence \((H_{n,k_\ell})_{\ell\geq1}\) converging locally uniformly in the spherical metric on \(B\) to a limit \(h\). 
    Since the full sequence tends to \(\infty\) on \(B\setminus\jul_n\), by the previous part, we must therefore have \(h\equiv\infty\) on this open set and hence on all of \(B\).
    This contradicts \(|H_{n,k_\ell}(u)|\leq1\) for any \(u\in B\cap\jul_n\subseteq \jul_n\).
    We conclude that the Julia set of \((H_{n,k})_{k\geq1}\) is \(\jul_n\).

    Now as \(\mu_n\) is fixed at a fixed phase \(n\), and as \(\widetilde F_{n,k}=H_{n,k}\circ\mu_n\) we get that $\widetilde\jul_n=\mu_n^{-1}(\jul_n)=A_{n-1}^{-1}\jul_n$.
    In particular, \(\widetilde\jul_n\) is a compact planar Cantor set, and dilation gives $\area(\widetilde\jul_n)=A_{n-1}^{-2}\area(\jul_n)>0$.
    Since \(A_0=1\), we have \(\widetilde\jul_1=\jul_1\) and \(\area(\widetilde\jul_1) =a\).

    Finally, fix $n\in\mbN$ and fix \(z\in\widetilde\jul_n\) and \(r>0\). 
    Then \(A_{n-1}z\in\jul_n\). 
    By \Cref{thm:survivor}, for all sufficiently large \(k\),
    \[
        F_{n,k}\left(B(A_{n-1}z,A_{n-1}r)\right)
        \supseteq\overline{\mbD}\supseteq\jul_{n+k}.
    \]
    Since \(F_{n,k}^{-1}(\jul_{n+k})=\jul_n\), every point of \(\jul_{n+k}\) has a preimage in \(\jul_n\cap B(A_{n-1}z,A_{n-1}r)\).
    The set \(\jul_{n+k}\) has positive area, and the polynomial  \(F_{n,k}\) maps sets of area zero to sets of area zero.
    Therefore \(\area(\jul_n\cap B(A_{n-1}z,A_{n-1}r))>0\). 
    Scaling back gives
    \[
        \area\bigl(\widetilde\jul_n\cap B(z,r)\bigr)
        =A_{n-1}^{-2}\area\bigl(\jul_n\cap B(A_{n-1}z,A_{n-1}r)\bigr)>0.
    \]
    This completes the proof. 
\end{proof}

%%%%%%%%%%%%%%%%%%%%%%%%%%%%%%%%%%%%%%%%%%%%%%%%%%%%%%%%%%%%%%%%%%%%%%%%%%%%%%%%%%%%%%%%%%%%%%%%%%%%
%%%%%%%%%%%%%%%%%%%%%%%%%%%%%%%%%%%%%%%%%%%%%%%%%%%%%%%%%%%%%%%%%%%%%%%%%%%%%%%%%%%%%%%%%%%%%%%%%%%%

\section{Quadratic sequences under affine changes of coordinates.}
\label{sec:quadratic}

In this section, we start with a fixed sequence of quadratic polynomials and change coordinates at each step by an affine map \(M_m(z)=m(z-a_m)\). 
We first prove \Cref{thm:quadratic-compact-pullback} for an arbitrary nonempty compact set \(K\).
We then use it to prove the three parts of \Cref{thm:main-quadratic} and to give the further example in \Cref{sec:quadratic-mixed}.

Let \(K\subset\mathbb C\) be a nonempty compact set, and choose a bounded sequence \((a_m)_{m\geq1}\) in \(\mbQ(i)\) whose set of subsequential limits (set of cluster points) is exactly \(K\). 
As observed in the introduction, we must necessarily have that \(\dist(a_m,K)\to0\).
Set \(M_0 := \mathrm{id}\) and \(M_m(z) := m(z-a_m)\) for \(m\geq1\). 
Then \(M_m\to\infty\) locally uniformly on \(\mathbb C\setminus K\). 
More precisely, whenever \(a_{m_j}\to a\), we have \(M_{m_j}\to\infty\) locally uniformly on \(\mathbb C\setminus\{a\}\).
Indeed, let \(E\subset\mbC\setminus\{a\}\) be a nonempty compact set and put \(\delta := \operatorname{dist}(E,a)>0\). 
Since \(a_{m_j}\to a\), for all sufficiently large \(j\) and every \(z\in E\) we have
\[
    |M_{m_j}(z)|
    =m_j|z-a_{m_j}|
    \geq m_j\bigl(\delta-|a_{m_j}-a|\bigr)
    \geq \tfrac12 m_j\delta.
\]
The last bound tends to $\infty$ uniformly on \(E\).

Let \(q_m(z) := z+\epsilon_m z^2\), with $\epsilon_m$ to be determined below. 
As the proof of \Cref{thm:quadratic-compact-pullback} below shall show, it suffices to choose positive rational numbers \(\epsilon_m\) such that
\[
    \sum_{m\geq1}\epsilon_m A^{2^m}<\infty
    \quad\text{for every } \ A>1,
    \qquad\text{ and  with}\qquad
    \sum_{m\geq1}\epsilon_m<\frac1{32}.
\]
For the explicit formulas below, we therefore make the choice \(\epsilon_m := 2^{-2^{m^2+4}}\).
Thus
\begin{equation}
\label{eq:q-sequence}
      q_m(z)=z+2^{-2^{m^2+4}} z^2   
\end{equation}
With this choice, for the sequence \((a_m)_{m\ge 1}\), define \(p_m := M_m\circ q_m\circ M_{m-1}^{-1}\),
where \(M_0=\id\) and \(M_m(z)=m(z-a_m)\) for \(m\geq1\).
Each \(p_m\) is quadratic and has coefficients in \(\mbQ(i)\), since both affine changes of coordinates have coefficients in \(\mbQ(i)\).
For \(m\geq s\geq1\), we have that
\[
    p_m\circ\cdots\circ p_s
    =M_m\circ(q_m\circ\cdots\circ q_s)\circ M_{s-1}^{-1},
\]
since the intermediate affine maps cancel out. 
For background on phase-dependent affine conjugacies in polynomial iteration, see \cite{Comerford03}.
Phase-dependent affine changes of coordinates need not carry the original phase Julia sets onto those of the transformed sequence; see \cite[Example~2.3]{Mayer2013}.
Here the family $(M_m)_{m\ge 1}$ is not equicontinuous in the spherical metric.

Since Julia sets in the following construction contain \(\infty\), we use spherical area to state the area conclusion in neighborhoods of all their points. 
For a measurable set \(E\subseteq\CC\), its normalized spherical area is
\begin{equation}
\label{eq:spherical-measure}
     m_{\mathrm{sph}}(E)
     := \int_{E\cap\mbC}\frac{1}{\pi(1+|z|^2)^2}\,d\area(z),
\end{equation}
where \(\area\) denotes planar Lebesgue measure. 
The density is positive on \(\mbC\) and \(m_{\mathrm{sph}}(\{\infty\})=0\), so \(E\) has positive spherical area if and only if \(E\cap\mbC\) has positive planar area.
Moreover 
\[
    m_{\mathrm{sph}}(\mathbb C)
        =\int_0^{2\pi}\int_0^\infty
        \frac{r}{\pi(1+r^2)^2}\,dr\,d\theta
        \, =2\left[-\frac{1}{2(1+r^2)}\right]_0^\infty
        =1.
\]
Since \(m_{\mathrm{sph}}(\{\infty\})=0\), the full Riemann sphere also has normalized spherical area \(1\); see \cite[Chapter~19]{Bonk2010}. 

%%%%%%%%%%%%%%%%%%%%%%%%%%%%%%%%%%%%%%%%%%%%%%%%%%%%%%%%%%%%%%%%%%%%%%%%%%%%%%%%%%%%%%%%%%%%%%%%%%%%
%%%%%%%%%%%%%%%%%%%%%%%%%%%%%%%%%%%%%%%%%%%%%%%%%%%%%%%%%%%%%%%%%%%%%%%%%%%%%%%%%%%%%%%%%%%%%%%%%%%%

\subsection{Proof of Theorem~C}

Write \(\seq{p} := (p_m)_{m\geq1}\) and denote its phase-\(n\) Julia set by \(\jul_n^{\seq{p}}\). For \(n,k\geq1\), denote the \(k\)-fold compositions starting at phase \(n\) by
\[
    Q_{n,k} := q_{n+k-1}\circ\cdots\circ q_n,\qquad\text{ and }\qquad
    P_{n,k} := p_{n+k-1}\circ\cdots\circ p_n,
\]
with \(Q_{n,0} := \mathrm{id}_{\CC}\) and \(P_{n,0} := \mathrm{id}_{\CC}\). Fix a phase \(n\) and put
\[
    h_{n,k} := M_{n+k-1}\circ Q_{n,k}\qquad\text{ for each } \quad k\geq1.
\]
Denote by \(\jul_n^{\seq{h}}\) the set of points of \(\CC\) at which \((h_{n,k})_{k\geq1}\) fails to be normal. 
The intermediate affine maps cancel, giving \(P_{n,k}=h_{n,k}\circ M_{n-1}^{-1}\).
As in the proof of \Cref{prop:exterior}, \(M_{n-1}\) is fixed at phase \(n\), so
\[
    \jul_n^{\seq{p}}=M_{n-1}\bigl(\jul_n^{\seq{h}}\bigr).
\]
We write \(J_n := \jul_n^{\seq{p}}\) when the sequence is clear.

\QuadraticCompactPullback*
\begin{proof}
    Fix \(n\geq1\). 
    Let \(z\in\mbC\), let \(k\geq1\), and put \(w := Q_{n,k-1}(z)\).
    \begin{equation}
    \label{eq:Q-difference-estimate}
        Q_{n,k}(z)-Q_{n,k-1}(z) 
        =  q_{n+k-1}(w) - w
        =\epsilon_{n+k-1}\bigl(Q_{n,k-1}(z)\bigr)^2.
    \end{equation}
    To bound this increment, observe that \(0<\epsilon_m\leq1\) gives
    \[
        1+|q_m(w)|\leq 1+|w|+|w|^2\leq(1+|w|)^2,
    \]
    for all $w\in\mbC$. 
    Thus \(1+|Q_{n,j+1}(z)|\leq(1+|Q_{n,j}(z)|)^2\) for \(j\geq0\).
    Starting from \(Q_{n,0}(z)=z\), induction gives \(1+|Q_{n,j}(z)|\leq(1+|z|)^{2^j}\) for every \(j\geq0\).
    For \(|z|\leq R\), taking \(j=k-1\) gives \(\bigl|Q_{n,k-1}(z)\bigr|^2\leq(1+R)^{2\cdot2^{k-1}} =(1+R)^{2^k}\). 
    Hence \eqref{eq:Q-difference-estimate} yields
    \[
        \sup_{|z|\leq R}|Q_{n,k}(z)-Q_{n,k-1}(z)|
        \leq\epsilon_{n+k-1}(1+R)^{2^k}.
    \]
    Since \(k\leq n+k-1\), the sum of these bounds is at most \(\sum_{m\geq n}\epsilon_m(1+R)^{2^m}<\infty\) by the choice of \(\epsilon_m\). 
    By summing the successive differences, we see that  \((Q_{n,k})_{k\geq0}\) is uniformly Cauchy on \(\{|z|\leq R\}\).
    Since \(R>0\) was arbitrary, the sequence converges locally uniformly on \(\mbC\) to an \emph{entire} function \(U_n\), by Weierstrass's theorem; \cite[\S~1.4, p.~9]{Schiff1993}.
    Since each \(q_m\) fixes \(0\) and has derivative \(1\) there, the chain rule gives \(Q_{n,k}'(0)=1\) for every \(k\). 
    By the same Weierstrass theorem, \(U_n'(0)=\lim_{k\to\infty}Q_{n,k}'(0)=1\). 
    Hence \(U_n\) is nonconstant and so is \(H_n=U_n\circ M_{n-1}^{-1}\). 

    To identify \(\jul_n^{\seq{p}}\), recall that  \(\jul_n^{\seq{p}}=M_{n-1}(\jul_n^{\seq{h}})\).
    We first show that $\jul_n^{\seq{h}}\cap\mbC=U_n^{-1}(K)$.
    
    Since \(h_{n,k}=(n+k-1)(Q_{n,k}-a_{n+k-1})\) and \(\dist(a_{n+k-1},K)\to0\), the maps \(h_{n,k}\) tend
    locally uniformly to \(\infty\) on \(\mbC\setminus U_n^{-1}(K)\).
    To see this,  fix a nonempty compact set \(E\subset\mbC\setminus U_n^{-1}(K)\).
    The compact sets \(U_n(E)\) and \(K\) are disjoint, so \(\delta:=\operatorname{dist}(U_n(E),K)>0\).
    Since \(Q_{n,k}\to U_n\) locally uniformly, we have \(\sup_{w\in E}|Q_{n,k}(w)-U_n(w)|<\delta/4\) for all sufficiently large \(k\). 
    For such sufficiently large \(k\), we also have \(\dist(a_{n+k-1},K)<\delta/4\), so \(|U_n(w)-a_{n+k-1}|\geq3\delta/4\) for every \(w\in E\).
    \[
        |Q_{n,k}(w)-a_{n+k-1}|
        \geq |U_n(w)-a_{n+k-1}|-|Q_{n,k}(w)-U_n(w)|
        \geq \frac{\delta}{2}.
    \]
    Hence \(\inf_{w\in E}|h_{n,k}(w)|\geq(n+k-1)\delta/2\to\infty\), as $k\to\infty$, as claimed.
    This gives us that $ \jul_n^{\seq{h}}\cap\mbC\subseteq U_n^{-1}(K)$. 

    For the reverse inclusion, let \(w_0\in U_n^{-1}(K)\) and put \(a := U_n(w_0)\in K\). 
    We show that \((h_{n,k})_{k\geq1}\) is not normal at \(w_0\). 
    Since the tail \((a_m)_{m\geq n}\) has cluster set exactly \(K\), we may choose \(k_j\to\infty\) such that \(a_{n+k_j-1}\to a\), as $j \to \infty$.
    Then \(Q_{n,k_j}-a_{n+k_j-1}\to U_n-a\), as $j \to \infty$, locally uniformly. 
    For brevity, write \(m_j := n+k_j-1\), and put
    \[
        g_j(w) := Q_{n,k_j}(w)-a_{m_j},
        \qquad g(w) := U_n(w)-a.
    \]
    Then \(g_j\to g\) locally uniformly and \(h_{n,k_j}=m_jg_j\).
    The function \(g\) vanishes at \(w_0\) and is nonconstant (since $U_n$ is not constant).
    So we can find some \(r_0>0\) so that \(w_0\) is the only zero of \(g\) in \(\overline{B(w_0,r_0)}\).
    Fix \(0<r\leq r_0\). Since \(g\) has no zeros on
    \(\{|w-w_0|=r\}\), we have
    \[
        \delta_r:=\min_{|w-w_0|=r}|g(w)|>0.
    \]
    By uniform convergence on this circle, for all sufficiently large \(j\),
    \[
        |g_j(w)-g(w)|<\delta_r\leq |g(w)|,
    \]
    for $|w-w_0|=r$.
    By Rouch\'e's theorem, \(g_j\) and \(g\) have the same number of zeros in \(B(w_0,r)\), counted with multiplicity. 
    In particular,\(g_j\) also has a zero in \(B(w_0,r)\), for all sufficiently large $j$.

    Applying the preceding argument on disks shrinking to \(w_0\) and passing to a subsequence, we obtain points \(w_j\to w_0\) such that \(h_{n,k_j}(w_j)=m_jg_j(w_j)=0\). 
    On the punctured disk \(0<|w-w_0|<r_0\), the limit \(g\) has no zeros. 
    Hence \(g_j\to g\) locally uniformly and \(m_j\to\infty\) imply that \(h_{n,k_j}=m_jg_j\to\infty\) locally uniformly there.
    If the family were normal at \(w_0\), a further subsequence would converge locally uniformly in the spherical metric on a disk about \(w_0\). 
    Its limit would equal \(\infty\) on the punctured disk and hence at \(w_0\) by continuity. 
    This contradicts the zeros \(h_{n,k_j}(w_j)=0\) with \(w_j\to w_0\). Thus \(w_0\in\jul_n^{\seq{h}}\), proving the reverse inclusion and hence $\jul_n^{\seq{h}}\cap\mbC=U_n^{-1}(K)$.

    Since \(M_{n-1}\) is affine and \(H_n=U_n\circ M_{n-1}^{-1}\), the identity \(\jul_n^{\seq{p}}=M_{n-1}(\jul_n^{\seq{h}})\) gives
    \[
        \jul_n^{\seq{p}}\cap\mbC
        =M_{n-1}\bigl(U_n^{-1}(K)\bigr)
        =H_n^{-1}(K).
    \]
    It remains to show that \(\infty\in\jul_n^{\seq p}\).
    Since \(\jul_n^{\seq p}\cap\mbC=H_n^{-1}(K)\) and \(\jul_n^{\seq p}\) is closed in \(\CC\), it suffices to show that \(H_n^{-1}(K)\) is unbounded.
    As \(K\) is nonempty, it suffices in turn to find one \(a\in K\) for which \(H_n^{-1}(\{a\})\) is unbounded.
    We prove the stronger assertion: for every \(a\in\mbC\), the fiber \(U_n^{-1}(\{a\})\), and hence its affine image \(H_n^{-1}(\{a\})\), is unbounded.
    
    Fix \(a\in\mbC\).
    Summing the increment estimate \eqref{eq:Q-difference-estimate} for the tail starting at \(m+1\) gives, for every \(R>0\),
    \[
        \sup_{|z|\leq R}|U_{m+1}(z)-z|
        \leq\sum_{\ell>m}\epsilon_\ell(1+R)^{2^\ell}
        \xrightarrow[m\to\infty]{} 0.
    \]
    On \(\partial B(a,1)\) we therefore have \(|U_{m+1}(z)-z|<1=|z-a|\), for all sufficiently large \(m\geq n\).
    By Rouch\'e's theorem there is \(w_m\in B(a,1)\) such that \(U_{m+1}(w_m)=a\).
    The two roots of \(q_m(v)=w_m\), counted with multiplicity, sum to \(-1/\epsilon_m\), so one root \(v_m\) satisfies \(|v_m|\geq1/(2\epsilon_m)\).
    Choose \(z_m\) with \(Q_{n,m-n}(z_m)=v_m\).
    Since \(Q_{n,m-n}\to U_n\) locally uniformly and \(|v_m|\to\infty\), the points \(z_m\) eventually leave every compact subset of \(\mbC\).
    Passing to the limit in the composition identity gives \(U_n=U_{m+1}\circ q_m\circ Q_{n,m-n}\), so \(U_n(z_m)=a\).
    Thus \(H_n^{-1}(\{a\})=M_{n-1}(U_n^{-1}(\{a\}))\) is unbounded, as is \(H_n^{-1}(K)\) because \(K\) is nonempty.
    Consequently, \(\infty\) is an accumulation point of \(\jul_n^{\seq{p}}\cap\mbC\). 
    Since \(\jul_n^{\seq{p}}\) is closed in \(\CC\), we conclude that
    \[
        \jul_n^{\seq{p}}=H_n^{-1}(K)\cup\{\infty\}.
    \]

    Let \(V\subseteq\CC\) be open and meet \(\jul_n^{\seq p}\).
    Put \(E := V\cap H_n^{-1}(K)\) and \(W := H_n(V\cap\mbC)\).
    Then \(E\) is the set of finite Julia points in \(V\), and \(H_n(E)=W\cap K\).
    We compare \(E\) with its image under \(H_n\).

    Let \(z\in E\) satisfy \(H_n'(z)\ne0\).
    By the holomorphic inverse function theorem \cite[Theorem~10.30]{Rudin1987}, choose a local inverse branch \(\phi\) of \(H_n\) on a disk \(B\) about \(H_n(z)\) such that \(\phi(H_n(z))=z\).
    Shrinking \(B\) if necessary, we may assume that \(\phi(B)\subseteq V\) and that \(\phi\) is bi-Lipschitz on \(B\).
    Since \(E=V\cap H_n^{-1}(K)\), we have
    \[
        E\cap\phi(B)=\phi(B\cap K).
    \]
    Thus \(E\cap\phi(B)\) and \(B\cap K\) have the same Hausdorff dimension, and one has planar area zero if and only if the other does.

    Since \(\mbC\) has a countable basis, countably many such neighborhoods \(\phi(B)\) cover all points of \(E\) at which \(H_n'\ne0\).
    The remaining points of \(E\) lie in the zero set of \(H_n'\), which is at most countable because \(H_n'\) is entire and not identically zero.
    Their images are also at most countable.
    Taking the union of the pieces above and adjoining these exceptional points, we conclude that \(E\) and \(H_n(E)=W\cap K\) have the same Hausdorff dimension and are simultaneously of zero or positive planar area.

    On \(\mbC\), the spherical and Euclidean metrics give the same Hausdorff dimension, since they are locally bi-Lipschitz equivalent.
    Moreover, spherical area has positive density with respect to planar area by \eqref{eq:spherical-measure}.
    If \(\infty\in V\), then \(E\) is nonempty because \(H_n^{-1}(K)\) is unbounded, and \(V\cap\jul_n^{\seq p}=E\cup\{\infty\}\).
    Adding this one point changes neither the Hausdorff dimension nor the spherical area.
    This proves both assertions.

    Finally, suppose that a finite point \(z\) lies in the interior of \(\jul_n^{\seq p}\).
    Choose an open disk \(D\) about \(z\) contained in \(\jul_n^{\seq p}\cap\mbC\).
    By the open mapping theorem, \(H_n(D)\) is an open neighborhood of \(H_n(z)\), and it lies in \(K\) because \(D\subseteq H_n^{-1}(K)\).
    Thus \(H_n(z)\in\operatorname{int}(K)\).
    Conversely, if \(H_n(z)\in\operatorname{int}(K)\), then \(H_n^{-1}(\operatorname{int}K)\) is an open neighborhood of \(z\) contained in \(\jul_n^{\seq p}\).
    This proves the interior equivalence.
\end{proof}

%%%%%%%%%%%%%%%%%%%%%%%%%%%%%%%%%%%%%%%%%%%%%%%%%%%%%%%%%%%%%%%%%%%%%%%%%%%%%%%%%%%%%%%%%%%%%%%%%%%%
%%%%%%%%%%%%%%%%%%%%%%%%%%%%%%%%%%%%%%%%%%%%%%%%%%%%%%%%%%%%%%%%%%%%%%%%%%%%%%%%%%%%%%%%%%%%%%%%%%%%

\subsection{Positive-area Cantor Julia set construction}

Let \(\mcC\subset[0,1]\) be the Smith--Volterra--Cantor (fat Cantor) set. 
At step \(j\geq1\), remove the open middle interval of length \(4^{-j}\) from each interval retained at step \(j-1\). 
Each of the two remaining intervals has length \(l_j := 2^{-j-1}+2^{-2j-1}\), and their left endpoints differ by \(c_j := l_j+4^{-j}=2^{-j-1}+3\,2^{-2j-1}\).
Writing \(m_1\) for one-dimensional Lebesgue measure, we have
\[
    m_1(\mcC)=\lim_{j\to\infty}2^jl_j=\frac12.
\]
Thus \(K := \mcC+i\mcC\) has planar area \(1/4\).

For \(m\geq1\), write \(m-1=\sum_{\nu\geq0}b_\nu(m)2^\nu\), with \(b_\nu(m)\in\{0,1\}\), and put
\[
    a_m := \sum_{j\geq1}\bigl(b_{2j-2}(m)+i b_{2j-1}(m)\bigr)c_j.
\]
Fix \(z\in K\) and \(j\geq1\). The first \(2j\) binary digits of \(m-1\), starting with the least significant digit, determine the square of side \(l_j\) containing \(a_m\) at step \(j\). 
Choosing these digits so that the square contains \(z\) gives \(|a_m-z|\leq\sqrt2\,l_j\).
Since the higher digits are free, this holds for arbitrarily large \(m\).
As \(l_j\to0\) and \(a_m\in K\) for every \(m\), each tail of \((a_m)\) has cluster set \(K\).

Take \(M_m\), \(q_m\) and \(p_m\) as defined at the beginning of \Cref{sec:quadratic}.
The sum defining \(a_m\) has only finitely many nonzero terms, and each \(c_j\) is rational. 
Thus \(a_m\in\mbQ(i)\), and since \(\epsilon_m\in\mbQ\), we have \(p_m\in\mbQ(i)[z]\).
For \(m\geq2\), substituting \(M_{m-1}^{-1}(z)=z/(m-1)+a_{m-1}\) into \(p_m=M_m\circ q_m\circ M_{m-1}^{-1}\) gives
\[
    p_m(z)=m\left[
        \frac{z}{m-1}+a_{m-1}-a_m+
        \epsilon_m\left(\frac{z}{m-1}+a_{m-1}\right)^2
    \right].
\]
For example, \(a_1=0\) and \(a_2=c_1=5/8\), so the first two maps for the explicit choice \(\epsilon_m=2^{-2^{m^2+4}}\) are
\[
    p_1(z)=z+2^{-32}z^2,\qquad p_2(z)=2z+2^{-255}z^2-\frac54.
\]

\begin{prop}
\label{prop:quadratic-fat}
    For the explicit quadratic sequence \(\seq p\) defined above, each phase Julia set \(\jul_n^{\seq p}\) is a Cantor set containing \(\infty\).
    More precisely, every open set \(V\subseteq\CC\) meeting \(\jul_n^{\seq p}\) satisfies
    \[
        m_{\mathrm{sph}}\bigl(V\cap\jul_n^{\seq p}\bigr)>0.
    \]
\end{prop}
\begin{proof}
    Fix a phase \(n\).
    Since the sequence \((a_m)\) has cluster set \(K=\mcC+i\mcC\), \Cref{thm:quadratic-compact-pullback} gives
    \[
        \jul_n^{\seq p}=H_n^{-1}(K)\cup\{\infty\}.
    \]

    At step \(j\) in the construction of \(K\), the \(4^j\) retained squares meet \(K\) in translates of the same set.
    Each of these sets therefore has area \(4^{-j}\area(K)>0\).
    Since the diameters of the squares tend to zero, every neighborhood of a point of \(K\) meets \(K\) in positive planar area.
    Now let \(V\subseteq\CC\) be open and meet \(\jul_n^{\seq p}\).
    Since \(H_n^{-1}(K)\) is unbounded, \(V\) contains a finite Julia point even if it is a neighborhood of \(\infty\).
    Thus \(H_n(V\cap\mbC)\) is open and meets \(K\).
    The local area equivalence in \Cref{thm:quadratic-compact-pullback} gives \(m_{\mathrm{sph}}(V\cap\jul_n^{\seq p})>0\).
    In particular, \(\jul_n^{\seq p}\) has no isolated points.

    It remains to prove total disconnectedness.
    Let \(C\) be a nonempty connected subset of \(H_n^{-1}(K)\).
    Since \(H_n(C)\) is connected and \(K\) is totally disconnected, \(H_n(C)=\{a\}\) for some \(a\in K\).
    The fiber \(H_n^{-1}(\{a\})\) is discrete because \(H_n\) is nonconstant and entire, so \(C\) is a singleton.
    As a closed subset of \(\mbC\), the totally disconnected space \(H_n^{-1}(K)\) is locally compact and hence has a basis of compact open sets \cite{willard_72}.
    Each such set is clopen in \(H_n^{-1}(K)\cup\{\infty\}\), so these sets separate any two distinct Julia points, including a finite point and \(\infty\).
    Thus \(\jul_n^{\seq p}\) is compact, perfect, and totally disconnected, and hence is a Cantor set.
\end{proof}

For the positive-area monomial examples in \cite[\S~5.2]{monomial}, the radial directions give an invariant sequence of non-trivial measurable line fields. 
Comerford~\cite{Comerford2013} constructed a bounded sequence of quadratic polynomials whose phase Julia sets also carry such fields. 
The following corollary gives a quadratic example outside Comerford's bounded class, since the leading coefficients \(m\epsilon_m/(m-1)^2\) of our maps \(p_m\) tend to zero. 

Recall that a measurable line field on a measurable set \(E\subseteq\CC\) assigns an unoriented tangent line to each point of \(E\). 
In a local holomorphic coordinate, the line spanned by a nonzero vector \(v\) is represented by the Beltrami differential \(\mu(z)\,d\bar z/dz\), with coefficient \(\mu(z) := v/\overline v\). 
This coefficient has modulus \(1\) on \(E\) and is extended by \(0\) outside \(E\). A sequence of such fields is invariant if \(p_s^*\mu_{s+1}=\mu_s\) almost everywhere for every \(s\); see \cite[\S~1.2, Definition~1.1]{Comerford2013}.

\begin{lem}
\label{cor:quadratic-line-fields}
    Let \(K\) and \((a_m)\) satisfy the hypotheses of \Cref{thm:quadratic-compact-pullback}, and suppose that every open disk meeting \(K\) meets it in positive planar area.
    The quadratic sequence \(\seq{p}\) admits an invariant sequence \((\mu_s)_{s\geq1}\) of non-trivial measurable line fields.
    At each phase \(s\), the essential support of \(\mu_s\) is \(\jul_s^{\seq{p}}\), and \(\mu_s\) admits a representative assigning a line at all but countably many points of \(\jul_s^{\seq{p}}\).
\end{lem}
\begin{proof}
    For each \(s\), let
    \[
        E_s := \{z\in\mbC:H_s(z)\in K,\ H_s'(z)\ne0\}.
    \]
    At each \(w\in K\), take the horizontal tangent line, spanned by \(1\).
    For \(z\in E_s\), its pullback under \(H_s\) is the line spanned by \(v_s(z) := 1/H_s'(z)\), since \(H_s'(z)v_s(z)=1\). 
    Its Beltrami coefficient is therefore
    \[
        \mu_s(z) := \frac{\, v_s(z)\,}{\overline{v_s(z)}}
                =\frac{\,\overline{H_s'(z)}\,}{H_s'(z)}.
    \]
    Set \(\mu_s := 0\) on \(\CC\setminus E_s\).
    The set \(E_s\) is Borel, and the Beltrami coefficient is measurable and has modulus \(1\) on \(E_s\). 
    By \Cref{thm:quadratic-compact-pullback}, \(\jul_s^{\seq{p}}=H_s^{-1}(K)\cup\{\infty\}\).
    Thus \(\jul_s^{\seq{p}}\setminus E_s\) consists of \(\infty\) and the finite points \(z\in\jul_s^{\seq{p}}\) at which \(H_s'(z)=0\).
    Since \(H_s\) is nonconstant and entire, this exceptional set is at most countable.
    
    To check invariance, pass to the limit, as $k\to\infty$, in \(Q_{s,k+1}=Q_{s+1,k}\circ q_s\). 
    This gives \(U_s=U_{s+1}\circ q_s\), and hence \(H_s=H_{s+1}\circ p_s\).
    Differentiating, we obtain
    \[
        H_s'(z)=H_{s+1}'(p_s(z))\,p_s'(z).
    \]
    At every finite \(z\) with \(p_s'(z)\ne0\), these two identities give  \(z\in E_s\) if and only if \(p_s(z)\in E_{s+1}\). 
    On \(E_s\),  the derivative identity yields
    \[
        \mu_{s+1}(p_s(z))
        \frac{\overline{p_s'(z)}}{\,p_s'(z)\,}=\mu_s(z).
    \]
    Outside \(E_s\), both sides vanish wherever \(p_s'(z)\ne0\).
    Thus \(p_s^*\mu_{s+1}=\mu_s\) almost everywhere on \(\CC\).
    
    By the hypothesis on \(K\) and the local area equivalence in \Cref{thm:quadratic-compact-pullback}, every open set meeting \(J_s\) meets it in positive spherical area.
    Since \(J_s\setminus E_s\) is at most countable, every such open set also meets \(E_s\) in positive spherical area.
    As \(|\mu_s|=1\) on \(E_s\) and \(\mu_s=0\) elsewhere, the field is non-trivial and its essential support is exactly \(J_s\).
\end{proof}

We now choose \(K\) to prescribe the spherical area at the first phase.

\begin{cor}
\label{cor:quadratic-prescribed-spherical-area}
    For every \(0<A<1\), there exists a Cantor set \(K_A\subset\mbC\) with the following property.
    For every bounded sequence \((a_m)_{m\geq1}\) in \(\mbQ(i)\) whose set of subsequential limits is \(K_A\), the quadratic sequence \((p_m)_{m\geq1}\) defined in \Cref{thm:quadratic-compact-pullback} satisfies \(m_{\mathrm{sph}}(\jul_1)=A\).
    At every phase \(s\), the set \(\jul_s\) is a Cantor set meeting every neighborhood of each of its points in positive spherical area, and the sequence admits invariant measurable line fields with essential support \(\jul_s\) at phase \(s\).
\end{cor}
\begin{proof}
    Since \(M_0=\id\), the map \(H_1=U_1\) is independent of the chosen sequence \((a_m)_m\) as long as its cluster points are exactly $K_A$. 
    For every Borel set \(E\subseteq\mbC\), put
    \[
        \nu(E):=m_{\mathrm{sph}}\bigl(U_1^{-1}(E)\bigr).
    \]
    Since \(U_1\) is entire, \(\nu(\mbC)=m_{\mathrm{sph}}(\mbC)=1\).
    As shown in the proof of \Cref{thm:quadratic-compact-pullback}, the preimage under \(U_1\) of a set of planar area zero has planar area zero and hence spherical area zero.
    Thus \(\nu\) is absolutely continuous with respect to \(\area\).

    Choose a closed square \(S\) with sides parallel to the coordinate axes and \(\nu(S)>A\).
    Write \(S=I+iJ\), where \(I\) and \(J\) are its projections onto the real and imaginary coordinates, respectively.
    Construct symmetric fat Cantor sets \(C_1\subset I\) and \(C_2\subset J\) by removing open middle intervals of arbitrarily small total length from each interval.
    Put \(C:=C_1+iC_2\).
    Every open disk meeting \(C\) meets it in positive planar area, similar to the argument carried out in \Cref{prop:quadratic-fat}.
    Since \(\area(S\setminus C)\) can be made arbitrarily small and \(\nu\) is finite and absolutely continuous with respect to \(\area\), we may choose \(C\) so that \(\nu(C)>A\).

    Since any line has \(\nu\)-measure zero,  the function
    \[
        t\mapsto \nu\bigl(C\cap\{\operatorname{Re}w<t\}\bigr)
    \]
    is continuous and takes all values between \(0\) and \(\nu(C)\).
    So there is some $t_0 \in \mbR$ such that $\nu\bigl(C\cap\{\operatorname{Re}w<t_0\}\bigr)=A$, and put
    \[
        K_A:=\overline{C\cap\{\operatorname{Re}w<t_0\}}.
    \]
    We still have that \(\nu(K_A)=A\).
    Every open disk meeting \(K_A\) contains a point \(w\in C\) with \(\operatorname{Re}w<t_0\).
    Then choose \(r>0\) small enough that
    \[
        B(w,r)\subseteq B\cap\{z\in\mbC:\operatorname{Re}z<t_0\}.
    \]
    Since \(w\in C\), the disk \(B(w,r)\) meets \(C\) in positive planar area.
    As \(B(w,r)\cap C\subseteq B\cap K_A\), it follows that \(\area(B\cap K_A)>0\).
    In particular, \(K_A\) has no isolated points, and as a compact subset of the totally disconnected set \(C\), it is itself a Cantor set.
    
    Now take any bounded sequence \((a_m)_{m\geq1}\) in \(\mbQ(i)\) whose set of subsequential limits is \(K_A\), and form the associated quadratic sequence \((p_m)_{m\geq1}\).
    Applying \Cref{thm:quadratic-compact-pullback} with \(K:=K_A\) gives
    \[
        \jul_s=H_s^{-1}(K_A)\cup\{\infty\}
        \qquad(s\geq1).
    \]
    At phase one, \(H_1=U_1\), so the definition of \(\nu\) and the fact that \(m_{\mathrm{sph}}(\{\infty\})=0\) give
    \[
        m_{\mathrm{sph}}(\jul_1)
        =m_{\mathrm{sph}}\bigl(U_1^{-1}(K_A)\bigr)
        =\nu(K_A)=A.
    \]
    
    Fix a phase \(s\) and an open set \(V\subseteq\CC\) meeting \(\jul_s\).
    The set \(W:=H_s(V\cap\mbC)\) is open and meets \(K_A\), also when \(V\) contains \(\infty\), since \(H_s^{-1}(K_A)\) is unbounded.
    By the local area property of \(K_A\) proved above, \(\area(W\cap K_A)>0\).
    The local area equivalence in \Cref{thm:quadratic-compact-pullback} therefore gives \(m_{\mathrm{sph}}(V\cap\jul_s)>0\).
    In particular, \(\jul_s\) has no isolated points.
    Since \(K_A\) is totally disconnected, the compact-open basis argument in the proof of \Cref{prop:quadratic-fat} shows that \(\jul_s\) is a Cantor set.
    The assertion about invariant line fields follows from \Cref{cor:quadratic-line-fields}.
\end{proof}

%%%%%%%%%%%%%%%%%%%%%%%%%%%%%%%%%%%%%%%%%%%%%%%%%%%%%%%%%%%%%%%%%%%%%%%%%%%%%%%%%%%%%%%%%%%%%%%%%%%%
%%%%%%%%%%%%%%%%%%%%%%%%%%%%%%%%%%%%%%%%%%%%%%%%%%%%%%%%%%%%%%%%%%%%%%%%%%%%%%%%%%%%%%%%%%%%%%%%%%%%

\subsection{Julia sets with countably infinitely many components}
\label{rem:countable-julia-components}

    Keep the quadratic polynomials \(q_m\) from \eqref{eq:q-sequence} and take
    \[
        D_j:=B(2^{-j},2^{-j-3}),
        \qquad\text{ and }\qquad 
        K:=\{0\}\cup\bigcup_{j\geq1}\overline{D_j}.
    \]
    Take any bounded sequence \((a_m)\) in \(\mbQ(i)\) with cluster set \(K\), and define \(M_m\) and \(p_m\) as in \Cref{thm:quadratic-compact-pullback}.

    We show that, at every phase \(s\), the Julia set \(\jul_s^{\seq{p}}\) has countably infinitely many connected components.

    Fix \(s\geq1\). By \Cref{thm:quadratic-compact-pullback}, \(\jul_s^{\seq{p}}=H_s^{-1}(K)\cup\{\infty\}\), where \(H_s=U_s\circ M_{s-1}^{-1}\).

    We first exhibit that each $\jul_s^{\seq{p}}$ has infinitely many distinct connected components. 
    Put \(z_s := M_{s-1}(0)\). Since \(U_s(0)=0\) and \(U_s'(0)=1\), we have $H_s(z_s)=0$ and 
    \[
        H_s'(z_s)=U_s'(0)(M_{s-1}^{-1})'(z_s)\ne0.
    \]
    By the holomorphic inverse function theorem \cite[Theorem~10.30]{Rudin1987}, \(H_s\) has a holomorphic inverse branch \(\psi_s\) on some disk \(B(0,r_s)\), with \(\psi_s(0)=z_s\).
    For all sufficiently large \(j\), we have \(\overline{D_j}\subset B(0,r_s)\), since \(2^{-j}+2^{-j-3}\to0\).
    For each such \(j\), the set \(\psi_s(\overline{D_j})\) is compact and connected.
    Since \(\overline{D_j}\) is separated from \(K\setminus \overline{D_j}\), there is an open set \(O_j\subset B(0,r_s)\) such that \(O_j\cap K=\overline{D_j}\). 
    Hence $\psi_s(\overline{D_j})=\psi_s(O_j)\cap\jul_s^{\seq{p}}$.
    The set \(\psi_s(O_j)\) is open, so \(\psi_s(\overline{D_j})\) is relatively open in \(\jul_s^{\seq{p}}\). 
    It is also relatively closed by compactness. 
    Thus \(\psi_s(\overline{D_j})\) is a connected component of \(\jul_s^{\seq{p}}\). 
    These components are distinct because the closed disks \(\overline{D_j}\) are pairwise disjoint and \(\psi_s\) is injective.

    For the upper bound, put
    \[
        A_s := H_s^{-1}(\{0\})
            \cup\{z\in\mbC:H_s'(z)=0\}\cup\{\infty\}.
    \]
    Since \(H_s\) is nonconstant and entire, \(A_s\) is at most countable.
    Let \(C\) be a component of \(\jul_s^{\seq{p}}\) containing a point \(z\notin A_s\). Then \(H_s(z)\in {\overline{D_j}}\) for some \(j\), and \(H_s'(z)\ne0\). 
    Choose a local inverse branch \(\phi\) of \(H_s\) on a disk \(B\) about \(H_s(z)\) such that \(\phi(H_s(z))=z\).
    The set \(B\cap {\overline{D_j}}\) is convex, so \(\phi(B\cap {\overline{D_j}})\) is a connected subset of \(\jul_s^{\seq{p}}\) containing \(z\). 
    It is therefore contained in \(C\). 
    Moreover, it contains the nonempty open set \(\phi(B\cap{D_j})\). Thus \(C\) has nonempty interior.
    There are at most countably many components with nonempty interior, since their interiors are pairwise disjoint and each contains a point of \(\mbQ(i)\). 
    All remaining components are contained in the countable set \(A_s\). Hence \(\jul_s^{\seq{p}}\) has countably infinitely many connected components.

    Finally, let \(V\subseteq\CC\) be an open set meeting \(\jul_s^{\seq{p}}\). 
    Then \(V\) contains a finite Julia point \(z\): if \(\infty\in V\), this follows from the unboundedness of \(H_s^{-1}(K)\). 
    Choose an open disk \(D\) about \(z\) contained in \(V\cap\mbC\). By the open mapping theorem, \(H_s(D)\) is an open neighborhood of \(H_s(z)\in K\). 
    Since
    \[
        K=\overline{\bigcup_{j\geq1}{D_j}},
    \]
    the set \(H_s(D)\) meets \({D_j}\) for some \(j\).
    It follows that \(D\cap H_s^{-1}({D_j})\) is a nonempty open subset of \(V\cap\jul_s^{\seq{p}}\cap\mbC\), and hence contains an open disk. 
    In particular, $\area\bigl(V\cap\jul_s^{\seq{p}}\cap\mbC\bigr)>0$ and $ m_{\mathrm{sph}}\bigl(V\cap\jul_s^{\seq{p}}\bigr)>0$.
    The open-disk assertion also gives \(J_s=\overline{\operatorname{int}(J_s)}\) in \(\CC\).
    Since \(K=\overline{\bigcup_jD_j}\), every open disk meeting \(K\) meets it in positive planar area, so \Cref{cor:quadratic-line-fields} gives the invariant line fields required in part~(2) of \Cref{thm:main-quadratic}.

\begin{cor}
\label{cor:countable-prescribed-spherical-area}
    For every \(A\in(0,1)\), the set \(K\) can be chosen so that \(m_{\mathrm{sph}}(\jul_1)=A\) and every phase Julia set has countably infinitely many connected components and is the closure of its interior.
    The associated quadratic sequence in \(\mbQ(i)[z]\) admits an invariant sequence of measurable line fields whose essential support at phase \(s\) is exactly \(\jul_s\).
\end{cor}
\begin{proof}
    Write \(K_0:=\{0\}\cup\bigcup_{j\geq1}\overline{D_j}\) for the compact set used above.
    As in the proof of \Cref{cor:quadratic-prescribed-spherical-area}, the measure \(\nu(E)=m_{\mathrm{sph}}(U_1^{-1}(E))\) has mass \(1\) and is absolutely continuous with respect to planar area.
    Choose \(r>0\) so small that \(\nu(\overline{B(0,2r)})<\min\{A,1-A\}\).
    Since \(K_0\subset\mbD\), the set \(rK_0\) lies in \(B(0,r)\).
    For \(t\geq2r\), put
    \[
        K_t:=rK_0\cup\{w\in\mbC:2r\leq |w|\leq t\}.
    \]
    Every circle is \(\nu\)-null, so \(t\mapsto\nu(K_t)\) is continuous.
    Its value at \(t=2r\) is \(\nu(rK_0)<A\), and its limit as \(t\to\infty\) is at least \(1-\nu(B(0,2r))>A\).
    Fix \(t>2r\) with \(\nu(K_t)=A\) and apply \Cref{thm:quadratic-compact-pullback} with \(K=K_t\).
    Since \(H_1=U_1\), we have \(m_{\mathrm{sph}}(\jul_1)=\nu(K_t)=A\).

    Fix a phase \(s\).
    As in the preceding example, a local inverse branch of \(H_s\) near \(0\) takes the disks \(r\overline{D_j}\), for all sufficiently large \(j\), to distinct connected components of \(\jul_s\).
    Indeed, these disks accumulate at \(0\), while the added annulus lies outside \(B(0,2r)\).
    The same argument gives the upper bound on the number of components.
    The only additional case is a regular point \(z\) for which \(H_s(z)\) belongs to the closed annulus.
    In every sufficiently small neighborhood of \(H_s(z)\), the annulus contains a closed disk of positive radius containing \(H_s(z)\), tangent to the boundary circle if necessary.
    A local inverse branch through \(z\) takes this disk to a connected subset of \(\jul_s\) with nonempty interior.
    Thus, as above, every component with empty interior is contained in the countable set
    \[
        H_s^{-1}(\{0\})\cup\{z\in\mbC:H_s'(z)=0\}\cup\{\infty\},
    \]
    while the components with nonempty interior are at most countable because their interiors are pairwise disjoint.
    Hence \(\jul_s\) has countably infinitely many components.

    Finally, \(K_t=\overline{\operatorname{int}(K_t)}\), so the open mapping argument above gives \(\jul_s=\overline{\operatorname{int}(\jul_s)}\) in \(\CC\).
    Every open disk meeting \(K_t\) therefore meets it in positive planar area, and \Cref{cor:quadratic-line-fields} gives the asserted invariant line fields.
\end{proof}

%%%%%%%%%%%%%%%%%%%%%%%%%%%%%%%%%%%%%%%%%%%%%%%%%%%%%%%%%%%%%%%%%%%%%%%%%%%%%%%%%%%%%%%%%%%%%%%%%%%%
%%%%%%%%%%%%%%%%%%%%%%%%%%%%%%%%%%%%%%%%%%%%%%%%%%%%%%%%%%%%%%%%%%%%%%%%%%%%%%%%%%%%%%%%%%%%%%%%%%%

\subsection{Cantor Julia sets of prescribed Hausdorff dimension}

\begin{thm}
\label{thm:quadratic-prescribed-dimension}
    For every \(d\in[0,2]\), there exists a sequence \(\seq{p}=(p_m)_{m\geq1}\) of quadratic polynomials in \(\mbQ(i)[z]\) such that, at every phase \(s\), the Julia set \(\jul_s^{\seq{p}}\) is a Cantor set in \(\CC\) containing \(\infty\), with
    \[
        m_{\mathrm{sph}}\bigl(\jul_s^{\seq{p}}\bigr)=0,
        \qquad
        \dim_{\mathrm H}\jul_s^{\seq{p}}=d.
    \]
Moreover, \(\dim_{\mathrm H}(V\cap\jul_s^{\seq p})=d\) for every open set \(V\subseteq\CC\) meeting \(\jul_s^{\seq p}\).
\end{thm}
\begin{proof}
    Fix \(d\in[0,2]\).
    Choose an infinite set \(S\subseteq\mbN\) with infinite complement (in $\mbN$) such that 
    \[
        \dfrac{\#(S\cap\{1,\ldots,n\})}{n} \to \frac d 2 
    \]
    as $n\to\infty$.
    In particular, we may take 
    \[
        S = 
        \begin{cases}
             \{k^2:k\geq1\} & \text{ for  d = 0},\\
             \{j\geq1:\lfloor jd/2\rfloor>\lfloor(j-1)d/2\rfloor\} & \text{ for } 0 < d< 2, \\
             \mbN\setminus\{k^2:k\geq1\}  & \text{ for } d =2. 
        \end{cases}
    \]

    Now for such a fixed set $S$, let \(C_S\subseteq[0,1]\) consist of all numbers admitting a binary expansion in which the digits at positions outside \(S\) are zero, while the digits at positions in \(S\) may be chosen freely from \(\{0,1\}\).
    Thus
    \[
        C_S := \left\{
            \sum_{j\in S}e_j2^{-j} :
            e_j\in\{0,1\}
        \right\}.
    \]
    It can be shown that $(C_s)$ is a Cantor set.
    This can be seen by first proving that the map
    \[
        \pi:\{0,1\}^S\longrightarrow C_S,
        \qquad\text{ so that }\qquad
        (e_j)_{j\in S}\longmapsto\sum_{j\in S}e_j2^{-j},
    \]
    is a homeomorphism. We omit these details.

    Now take 
    \[
        K:=C_S+iC_S
          =\{x+iy:x,y\in C_S\}.
    \]
    The set \(K\) is also a Cantor set. 
    It lies in the unit square \([0,1]+i[0,1]\) and contains \(0\), corresponding to choosing every digit in both coordinates to be zero.
    From \cite[Lemma~1]{WEI_2016}, we have that 
    \[
        \dim_{\mathrm H}K
        =\dim_{\mathrm H}(C_S\times C_S)
        =2\liminf_{n\to\infty}
          \frac{\#(S\cap\{1,\ldots,n\})}{n}
        =d.
    \]
    Moreover, the choices of the first \(n\) digits in both coordinates give a cover of \(K\) by \(4^{h_n}\) dyadic squares of side length \(2^{-n}\).
    Consequently,
    \[
        \area(K)
        \leq4^{h_n}(2^{-n})^2
        =4^{h_n-n}
        \longrightarrow0,
    \]
    so \(\area(K)=0\).

    We next verify that \(\dim_{\mathrm H}(U\cap K)=d\), for every open set \(U\subseteq\mbC\) that intersects \(K\).
    Let \(K_N\) be the subset of \(K\) obtained by fixing the first \(N\) digits in both coordinates to be zero.
    Every subset obtained by fixing the first \(N\) digits to any values satisfying the digit restrictions is a translate of \(K_N\). 
    Thus \(K\) is a finite union of translates of \(K_N\), thus
    \[
        \dim_{\mathrm H}K_N=\dim_{\mathrm H}K=d.
    \]
    Fix \(z\in U\cap K\). 
    Choose binary expansions of its two coordinates satisfying the digit restrictions.
    The subset of \(K\) obtained by fixing the first \(N\) digits to agree with these expansions contains \(z\), is a translate of \(K_N\), and has diameter  at most \(\sqrt{2}\,2^{-N}\).
    For sufficiently large \(N\), this subset lies in \(U\). 
    Therefore
    \[
        d\leq\dim_{\mathrm H}(U\cap K)
          \leq\dim_{\mathrm H}K=d,
    \]
    proving the local dimension assertion of $K$.

    The points of \(K\) having only finitely many nonzero digits in both coordinates belong to \(\mbQ(i)\)  and are dense in \(K\), as follows by truncating the binary expansions. 
    Enumerate these points as \(b_1,b_2,\ldots\), and form a sequence \((a_m)_{m\geq1}\) by concatenating the finite lists
    \[
        (b_1),\qquad
        (b_1,b_2),\qquad
        (b_1,b_2,b_3),\qquad\ldots.
    \]
    Then \((a_m)\) is a bounded sequence in \(\mbQ(i)\), every tail is dense in \(K\), and its set of subsequential limits is exactly \(K\).
    Keep the quadratic polynomials \(q_m\) from  \eqref{eq:q-sequence}. Set \(M_0=\mathrm{id}\) and,  for \(m\geq1\), define
    \[
        M_m(z):=m(z-a_m),
        \qquad
        p_m:=M_m\circ q_m\circ M_{m-1}^{-1}.
    \]
    Since \(q_m\) has rational coefficients and th affine maps \(M_m\) and their inverses have coefficients in \(\mbQ(i)\), each \(p_m\) is quadratic and belongs to \(\mbQ(i)[z]\).

    Fix \(s\geq1\) and write \(J_s:=\jul_s^{\seq{p}}\).
    By \Cref{thm:quadratic-compact-pullback},
    \[
        J_s=H_s^{-1}(K)\cup\{\infty\},
    \]
    where \(H_s\) is nonconstant and entire, and every fiber of \(H_s\) is unbounded.
    The finite part \(H_s^{-1}(K)\) has no isolated points, by the open mapping theorem and the
    perfectness of \(K\). Since \(0\in K\) and \(H_s^{-1}(\{0\})\) is unbounded, \(\infty\) is  also a limit point of \(J_s\).

    Furthermore, \(H_s^{-1}(K)\) is totally disconnected.
    Indeed, any connected subset has a connected image in \(K\), so it lies in a single fiber of \(H_s\), which is discrete. 
    The compact-open basis argument in the proof of \Cref{prop:quadratic-fat} then shows that \(J_s\) is totally disconnected as well.
    Since \(J_s\) is compact and has no isolated points, it is a Cantor set.

    Let \(V\subseteq\CC\) be open and meet \(J_s\), and put
    \[
        W:=H_s(V\cap\mbC).
    \]
    The set \(W\) is open and meets \(K\). Indeed, this is immediate if \(V\) contains a finite point of
    \(J_s\); if \(V\) contains \(\infty\), it follows  from the unboundedness of \(H_s^{-1}(\{0\})\).
    Then \Cref{thm:quadratic-compact-pullback}, together with the local dimension property of \(K\), gives
    \[
        \dim_{\mathrm H}(V\cap J_s)
        =\dim_{\mathrm H}(W\cap K)
        =d.
    \]
    Taking \(V=\CC\) yields  \(\dim_{\mathrm H}J_s=d\).
    Finally, the area assertion of the same theorem, applied with \(V=\CC\), and \(\area(K)=0\) imply
    \[
        m_{\mathrm{sph}}(J_s)=0.
    \]
    Since \(s\geq1\) was arbitrary, all conclusions hold at every phase.
\end{proof}

%%%%%%%%%%%%%%%%%%%%%%%%%%%%%%%%%%%%%%%%%%%%%%%%%%%%%%%%%%%%%%%%%%%%%%%%%%%%%%%%%%%%%%%%%%%%%%%%%%%%
%%%%%%%%%%%%%%%%%%%%%%%%%%%%%%%%%%%%%%%%%%%%%%%%%%%%%%%%%%%%%%%%%%%%%%%%%%%%%%%%%%%%%%%%%%%%%%%%%%%

\subsection{Other local geometries in Julia sets}
\label{sec:quadratic-mixed}

For the further example announced in the introduction, take
\[
    K := \overline{\mbD}\cup\bigl(3+\mcC+i\mcC\bigr)\cup[6,7]
\]
and fix a phase \(s\).
By \Cref{thm:quadratic-compact-pullback}, \(H_s^{-1}(\mbD)\) is a nonempty open subset of \(J_s\), and we may choose \(x_s,y_s\) with \(H_s(x_s)=3\) and \(H_s(y_s)=13/2\).
For each sufficiently small disk \(D\) about \(x_s\), the open neighborhood \(H_s(D)\) of \(3\) misses \(\overline{\mbD}\) and \([6,7]\).
By the nested-square argument in \Cref{prop:quadratic-fat}, \(H_s(D)\cap K\) has positive planar area.
It contains no open disk, since it lies in a Cantor set.
The local area and interior assertions of \Cref{thm:quadratic-compact-pullback} show that \(D\cap J_s\) has positive planar area and empty interior.
For each sufficiently small disk \(D\) about \(y_s\), the open set \(H_s(D)\) misses the other two pieces of \(K\), so \(H_s(D)\cap K\) is contained in \([6,7]\) and contains a nondegenerate interval about \(13/2\).
Hence it has zero planar area and Hausdorff dimension one.
Both properties pass to \(D\cap J_s\) by the same theorem.
In this example, as in part~(2) of \Cref{thm:main-quadratic}, the nonempty interior of the phase Julia sets excludes a common forward invariant hyperbolic domain; see \cite[Satz~2.2.3]{Bue1995}.

\smallskip

The same theorem also gives phase Julia sets with isolated points.
For \(K=\{0\}\), take \(a_m=0\) for every \(m\).
Then \Cref{thm:quadratic-compact-pullback} gives \(\jul_s=H_s^{-1}(\{0\})\cup\{\infty\}\) at every phase \(s\).
Since \(H_s\) is nonconstant and entire and its fibers are unbounded, \(\jul_s\) is countably infinite, every finite Julia point is isolated, and \(\infty\) is its only accumulation point.
Alternatively, take \(K=\{0\}\cup\overline{B(3,1)}\) and choose a bounded sequence \((a_m)\) in \(\mbQ(i)\) with cluster set \(K\).
Every point of \(H_s^{-1}(\{0\})\) is again isolated in \(\jul_s\), whereas \(H_s^{-1}(B(3,1))\) is a nonempty open subset of \(\jul_s\).
Both subsets are unbounded, so at every phase \(\infty\) is approached by both isolated Julia points and interior Julia points.

\medskip

We finish this section with a proof of  \Cref{thm:main-quadratic}

\begin{proof}[Proof of \Cref{thm:main-quadratic}]
    Use the fixed quadratic sequence \((q_m)_{m\geq1}\) in \eqref{eq:q-sequence} throughout.
    For part~(1), \Cref{cor:quadratic-prescribed-spherical-area} gives the prescribed first phase spherical area, the Cantor and local positive-area properties at every phase, and the invariant line fields with the stated essential supports.
    For part~(2), \Cref{cor:countable-prescribed-spherical-area} gives the prescribed first phase area, the component count, the density of the interior, and the invariant line fields.
    For part~(3), \Cref{thm:quadratic-prescribed-dimension} gives Cantor Julia sets of zero spherical area with the prescribed Hausdorff dimension in every open set meeting them.

    In each case, take a bounded sequence \((a_m)_{m\geq1}\) in \(\mbQ(i)\) with the indicated cluster set \(K\).
    The formula \(p_m=M_m\circ q_m\circ M_{m-1}^{-1}\) gives quadratic polynomials in \(\mbQ(i)[z]\), and \Cref{thm:quadratic-compact-pullback} gives \(\infty\in\jul_s\) at every phase \(s\).
\end{proof}

%%%%%%%%%%%%%%%%%%%%%%%%%%%%%%%%%%%%%%%%%%%%%%%%%%%%%%%%%%%%%%%%%%%%%%%%%%%%%%%%%%%%%%%%%%%%%%%%%%%%
%%%%%%%%%%%%%%%%%%%%%%%%%%%%%%%%%%%%%%%%%%%%%%%%%%%%%%%%%%%%%%%%%%%%%%%%%%%%%%%%%%%%%%%%%%%%%%%%%%%%

\section*{Acknowledgements}
    The author was supported via the Danish e-infrastructure Consortium (DeiC) Grant No.~5260-00014B. 
    Parts of this work were also supported by Villum Fonden Grant No.~25452, Grant No.~60842, and QMATH Center of Excellence Grant No.~10059. 
    The author also thanks  Carsten Lunde Petersen and Dimitris Vardakis for discussions about non-autonomous Julia sets.

%%%%%%%%%%%%%%%%%%%%%%%%%%%%%%%%%%%%%%%%%%%%%%%%%%%%%%%%%%%%%%%%%%%%%%%%%%%%%%%%%%%%%%%%%%%%%%%%%%%%
%%%%%%%%%%%%%%%%%%%%%%%%%%%%%%%%%%%%%%%%%%%%%%%%%%%%%%%%%%%%%%%%%%%%%%%%%%%%%%%%%%%%%%%%%%%%%%%%%%%%

\printbibliography

%%%%%%%%%%%%%%%%%%%%%%%%%%%%%%%%%%%%%%%%%%%%%%%%%%%%%%%%%%%%%%%%%%%%%%%%%%%%%%%%%%%%%%%%%%%%%%%%%%%%
%%%%%%%%%%%%%%%%%%%%%%%%%%%%%%%%%%%%%%%%%%%%%%%%%%%%%%%%%%%%%%%%%%%%%%%%%%%%%%%%%%%%%%%%%%%%%%%%%%%%

\end{document}